# The behavior of sequences of solutions to the Vafa-Witten equations

Clifford Henry Taubes[†]

Department of Mathematics
Harvard University
Cambridge, MA 02138

chtaubes@math.harvard.edu

ABSTRACT: The Vafa-Witten equations on an oriented Riemannian 4-manifold are first order, non-linear equations for a pair of connection on a principle SO(3) bundle over the 4-manifold and a self-dual 2-form with values in the associated Lie algebra bundle. The main theorem in this paper characterizes in part the behavior of sequences of solutions to the Vafa-Witten equations which have no convergent subsequence. The paper proves that a renormalization of a subsequence of the self-dual 2-form components converges on the complement of a closed set with Hausdorff dimension at most 2, with the limit being a harmonic 2-form with values in a real line bundle.

[†]Supported in part by the National Science Foundation

## 1. Introduction

Let X denote a connected, smooth, oriented 4-dimensional manifold with a Riemannian metric. The metric's Hodge star operator (an endomorphism of $\wedge^*(T^*X)$) maps 2-forms to 2-forms with square 1. The respective +1 and -1 eigenspaces on $\wedge^2 T^*X$ are denoted by $\Lambda^+$ and $\Lambda^-$ in what follows; both are oriented 3-dimensional vector bundles. Let $P \to X$ denote a chosen principle SU(2) or SO(3) bundle: and let $\mathfrak{G} \to X$ denote the associated bundle whose fiber is the Lie algebra of SU(2). The Vafa Witten equations (see [VW]) are equations for a pair $(A, \hat{a})$ with A being a connection on P and with $\hat{a}$ being a section of $\Lambda^+ \otimes \mathfrak{G}$. These equations have the schematic form

$$d_A \hat{a} = 0 \quad and \quad F_A^+ = \tfrac{1}{2} [\hat{a}; \hat{a}],$$

(1.1)

where $d_A$ is the exterior covariant derivative, where $F_A^+$ is the self-dual part of the curvature of the connection A, and where [ ; ] is a certain quadratic map from $\Lambda^+ \otimes \mathfrak{G}$ which is defined directly. Here is a useful definition (for calculations): Fix an oriented, orthonormal frame for $\Lambda^+$ on any given small radius ball in X. Write the components of $\hat{a}$ with respect to this frame as $\{\hat{a}_a\}_{a=1,2,3}$. Then, the components of $[\hat{a}; \hat{a}]$ on the chosen ball are $\{\tfrac{1}{\sqrt{2}} \varepsilon^{abc} [\hat{a}_b, \hat{a}_c]\}_{a=1,2,3}$ with [ , ] denoting the Lie bracket on the Lie algebra and with $\varepsilon^{abc}$ denoting the components of the completely anti-symmetric tensor on $\mathbb{R}^3$ with $\varepsilon^{123} = 1$. Note that the summation convention is used here (and below): Repeated indices are summed implicitly unless noted to the contrary. (Different oriented, orthonormal frames for $\Lambda^+$ give the same $[\hat{a}; \hat{a}]$ because SO(3) matrices are orthogonal with determinant 1.)

Of interest in what follows are sequences of solutions to (1.1) for a given manifold X and bundle P. Of particular concern: What can be said about the convergence (or not) of some sequences? The answer is known (thanks to Ben Mares [Ma].) when there is an a priori bound on the $L^2$ norm on X of the $\Lambda^+ \otimes \mathfrak{G}$ components from the sequence. This is asserted in the upcoming Theorem 1.1. By way of terminology: The $L^2$ norm of a section of $\wedge^2 T^*X \otimes \mathfrak{G}$ is the integral over X of the square of the pointwise norm; and the latter is defined using an ad-invariant norm on the Lie algebra of SU(2) (about which more is said momentarily) and the Riemannian metric. Also, $\mathbb{N}$ is used in the theorem to denote the set of positive integers.

**Theorem 1.1** (Mares [Ma]): *Let $\{(A_n, \hat{a}_n)\}_{n \in \mathbb{N}}$ denote a sequence of solutions to (1.1) on X with a uniform bound on the $L^2$ norms of the elements in $\{\hat{a}_n\}_{n \in \mathbb{N}}$ and, if X is not compact, a uniform bound on the $L^2$ norms of the curvature tensors also. There is a principle SU(2) or SO(3) bundle $P' \to X$ (as the case may be), a pair $(A', \hat{a}')$ of solutions to the Vafa-Witten equations on X for the bundle $P'$, and a finite set of points in X (to be* *denoted by Ξ). There is also a subsequence* $\Lambda \subset \mathbb{N}$ *and a corresponding sequence of bundle isomorphisms* $\{g_n: P'|_{X-\Xi} \to P|_{X-\Xi}\}_{n \in \Lambda}$ *such that* $\{(g_n^*A_n, g_n^*\hat{a}_n)\}_{n \in \Lambda}$ *converges to* $(A', \hat{a}')$ *in the* $C^\infty$ *topology on compact subsets of* X–Ξ.

When $\hat{a} = 0$, then the Vafa-Witten equations say that the curvature of A is anti-self dual. Theorem 1.1 is therefore a generalization of the following corollary of Karen Uhlenbeck's theorems in [U1], [U2]: Let $\{A_n\}_{n \in \mathbb{N}}$ denote a sequence of connections on P with anti-self dual curvature 2-form and, when X is not compact, a uniform bound for the $L^2$ norm of the curvature tensors. There is a principle SU(2) or SO(3) bundle $P' \to X$, a connection $A'$ on $P'$ with anti-self dual curvature, and a locally finite set of points in X (to be denoted by Ξ). There is also a subsequence $\Lambda \subset \mathbb{N}$ and a corresponding sequence of bundle isomorphisms $\{g_n: P'|_{X-\Xi} \to P|_{X-\Xi}\}_{n \in \Lambda}$ such that $\{g_n^*A_n\}_{n \in \Lambda}$ converges to $A'$ in the $C^\infty$ topology on compact subsets of X–Ξ.

With regard to $L^2$ norms: If X is not compact, then $a$ need not have a finite $L^2$ norm. Even so, $a$ is always $C^\infty$; and as a consequence $a$ has finite $L^2$ norm on any open set in X with compact closure (for example, small radius balls). This being the case, Theorem 1.1 (and the upcoming Theorem 1.2) can be invoked on these pre-compact submanifolds of the ambient non-compact manifold.

Theorem 1.1 is not the end of the story because there are pairs (X, P) which have sequences of solutions to the Vafa-Witten equation with no a priori bounds on the $L^2$ norms of the $\Lambda^+ \otimes \mathfrak{G}$ components. See [BW] for cases where X is a compact, Kähler surface. Such sequences can also be found for the case when X is $M \times S^1$ with M being a Riemannian 3-manifold and the metric is the product metric. These sequence come from non-convergent sequences of flat $Sl(2; \mathbb{C})$ connections on M.

A change of notation puts the focus on the $L^2$ norm of the $\Lambda^+ \otimes \mathfrak{G}$ component of the solution: The equations in (1.1) (when $a$ is not identically zero and has finite $L^2$ norm) are equivalent to equations for a triple $(r, (A, a))$ with $r$ being a positive real number, A being a connection on P, and $a$ being a section of $\Lambda^+ \otimes \mathfrak{G}$:

- $d_A a = 0$.
- $F_A^+ = \frac{1}{2} r^2 [a; a]$ .
- $\int_X |a|^2 = 1$ .

(1.2)

The equivalence with (1.1) is obtained by taking $r$ to be the $L^2$ norm of $\hat{a}$ and taking $a$ to be $r^{-1}\hat{a}$. Of interest are sequences of solutions to (1.2) with the sequence of $r$ values having no convergent subsequences. The upcoming Theorem 1.2 talks about these sequences; it is the central observation of this paper.

The upcoming theorem introduces the notion of a real line bundle over a given open set in X. This is defined to be the associated $\mathbb{R}$ bundle to a $\mathbb{Z}/2$ principle bundle over the set (a double cover of the set). A real line bundle defined in this way has a canonical fiber metric and canonical compatible flat connection. There is also a canonical exterior derivative operator acting on differential forms with values in the real line bundle (denoted by d); and there is a canonical covariant derivative on real line bundle valued tensors. (The latter is defined with the metric's Levi-Civita covariant derivative.)

By way of more notation; a given measurable function on X is said to be an $L^2_{1;loc}$ function if the square of the function and the square of the norm of its exterior derivative are integrable on any open set in X with compact closure. If X is compact, there is no distinction between $L^2_{1;loc}$ and the Sobolev space of functions that are square integrable and have square integrable derivative. Keep in mind that the vector space $L^2_{1;loc}$ has a natural Fréchet space topology.

The final bit of notation concerns the curvature of connections on X: Supposing that A is a connection on X, the upcoming theorem uses $F_A^-$ to indicate the anti-self dual part of its curvature 2-form. (The full curvature of A is denoted in what follows by $F_A$.)

**Theorem 1.2**: *Let $\{(r_n, (A_n, a_n)\}_{n \in \mathbb{N}}$ denote a sequence of solutions to (1.2) with $\{r_n\}_{n \in \mathbb{N}}$ having no convergent subsequences. There exists a continuous and $L^2_{1;loc}$ function on X to be denoted by $|v_\ddagger|$ and a subsequence $\Lambda \subset \mathbb{N}$ such that $\{|a_n|\}_{n \in \Lambda}$ converges to $|v_\ddagger|$ in the $L^2_1$ and $C^0$ topologies on compact subsets of X. If $|v_\ddagger|$ is not identically zero (which is necessarily the case if X is compact), then the what follows is true.*

- *The zero locus of $|v_\ddagger|$ (which is denoted henceforth by Z) is a nowhere dense set with Hausdorff dimension at most 2.*
- *The function $|v_\ddagger|$ is smooth on the complement of Z and uniformly Hölder continuous along Z in compact subsets of X.*
- *There is a real line bundle $\mathcal{I} \to X-Z$ and a smooth, closed, $\mathcal{I}$-valued self-dual 2 form on X–Z (a section of $\Lambda^+ \otimes \mathcal{I}$) to be denoted by $v_\ddagger$ whose norm is the function $|v_\ddagger|$.*
- *There is a $\Lambda$-indexed sequence of isometric homomorphisms from $\mathcal{I}$ to $\mathfrak{G}|_{X-Z}$ to be denoted by $\{\sigma_n\}_{n \in \Lambda}$ that obeys the following:*
  a) *The sequence $\{a_n - v_\ddagger \sigma_n\}_{n \in \Lambda}$ converges to zero in the $C^0$ topology on compact subsets of X–Z. In addition, if $U \subset X-Z$ is any open set with compact closure, then*
     i)   $\lim_{n \in \Lambda, n \to \infty} \int_U |\nabla_{A_n} \sigma_n|^2 = 0.$

ii) $\lim_{n\in\Lambda, n\to\infty} \int_U |\nabla_{A_n}(a_n - \nu_\ddagger \sigma_n)|^2 = 0$.

iii) $\lim_{n\in\Lambda, n\to\infty} \int_U |F_{A_n}^+|^2 = 0$.

b) *If* X *is compact, or if there is an upper bound for the sequence* $\{\frac{1}{r_n^2}\int_X |F_{A_n}^-|^2\}_{n\in\Lambda}$
*(which is the case if* X *is compact), and supposing that* $U\subset X-Z$ *is an open set with compact closure, then*

i) $\lim_{n\in\Lambda, n\to\infty} \int_U |\nabla_{A_n}\nabla_{A_n}\sigma_n|^2 = 0$.

ii) $\lim_{n\in\Lambda, n\to\infty} \int_U |\nabla_{A_n}\nabla_{A_n}(a_n - \nu_\ddagger \sigma_n)|^2 = 0$.

In the non-compact setting, it is entirely possible that $\nu_\ddagger$ is everywhere zero (its $L^2$ norm on X at most equal to 1.) In the compact setting, its $L^2$ norm is equal to 1 so $\nu_\ddagger$ is interesting. In any event, if $\nu_\ddagger$ is not identically zero, then the data $(Z, \mathcal{I}, \nu_\ddagger)$ from Theorem 1.2 defines what [T1] calls a $\mathbb{Z}/2$ harmonic spinor. Therefore, the theorems in [T1] say more about the structure of Z and the behavior of $\nu_\ddagger$ near Z.

What follows next are some additional remarks about Theorem 1.2.

REMARK 1: Yuuji Tanaka [Ta] proved a theorem to the effect that if the sequence of solutions $\{(A_n, \hat{a}_n)\}_{n\in\mathbb{N}}$ has uniformly small $L^2$ curvature norm on balls of fixed, small radius (the bound being a universal, small positive number), then either there is a convergent subsequence (after action of Aut(P)), or the sequence $\{A_n\}_{n\in\mathbb{N}}$ of connections has a subsequence that converges (after action of Aut(P)) in the $L^2_1$ weak topology to a reducible connection on P with self-dual curvature. (The weak $L^2_1$ convergence to some $L^2_1$ connection on P is directly a consequence of Karen Uhlenbeck's theorems in [U1].)

REMARK 2: Theorem 1.2 makes no claim to the effect that the sequence of connections $\{A_n\}_{n\in\Lambda}$ has a convergent subsequence on compact subsets X–Z even allowing for the action on each term of some automorphism of P. Nor is there a claim (which would be a weaker claim) that the sequence $\{\sigma_n\}_{n\in\Lambda}$ converges on compact subsets of X–Z after the term by term action of automorphisms of $\mathfrak{G}$. This is in stark contrast to what is said by the theorems in [T2] and [T3] which consider equations that are formally very similar to (1.2). More is said about this state of affairs in the upcoming Section 1a. In any event, there is a $\Lambda$–indexed sequence of automorphisms of P on any given ball in X–Z whose n'th member identifies $\sigma_n$ with a fixed section of $\mathfrak{G}$. Also, given the assumptions for Part c) of Theorem 1.2, then

$$\lim_{n \in \Lambda, n \to \infty} \int_U \|[F_{A_n}, \sigma_n]\|^2 = 0.$$

(1.3)

Thus, at least part of $A_n$'s curvature is small on X–Z when $n \in \Lambda$ is large. (Section 5d of this paper says slightly more about the curvature of the large $n \in \Lambda$ versions of $A_n$.)

REMARK 3: Adding terms to the right hand side of the curvature equation in (1.1) (and correspondingly to the equation in the second bullet of (1.2)) will not change the conclusions of Theorems 1.1 or 1.2 if these terms are relatively mild. For example, the Theorems 1.1 and 1.2 hold when (1.1) is changed to read

$$d_A \hat{a} = 0 \quad and \quad F_A^+ = \tfrac{1}{2}[\hat{a}; \hat{a}] + m\hat{a}$$

(1.4)

with $m$ being constant. The corresponding version of (1.2) reads

- $d_A a = 0.$
- $F_A^+ = \tfrac{1}{2} r^2 [a; a] + mra$.
- $\int_X |a|^2 = 1$.

(1.5)

The proof that Theorems 1.1 and 1.2 hold with (A, a) obeying (1.5) is almost the same as that with $m = 0$ (and isomorphic with regards to the sequence of steps). For this reason, only the $m = 0$ case is considered in what follows (except for some brief comments in Section 7a).

REMARK 4: Vafa and Witten consider a generalization of (1.1) (or (1.4)) of different sort: This generalization is a sysetm of equations for a triple $(A, (\hat{a}, \hat{\phi}))$ with A and $\hat{a}$ as before, and with $\hat{\phi}$ being a section of the bundle $\mathfrak{G}$. The equations in (1.1) are replaced by the following:

$$d_A \hat{a} + *d_A \hat{\phi} = 0 \quad and \quad F_A^+ = \tfrac{1}{2}[\hat{a};\hat{a}] + \tfrac{1}{\sqrt{2}}[\hat{\phi}, \hat{a}].$$

(1.6)

Adding $\hat{\phi}$ is not so interesting if X is compact because if X is compact, then either $\hat{\phi} = 0$, or $\hat{\phi}$ is A-covariantly constant and $[\hat{\phi}, \hat{a}] = 0$. In the former case, the equations are those in (1.2); and in the latter case, the connection A is reducible with anti-self dual curvature and $a$ can be written as $v_+ \hat{\phi}$ with $v_+$ being a self-dual harmonic 2-form. One can also add an $m\hat{a}$ term to the right hand side of (1.6) (mimicking (1.4)). When $m \neq 0$, it is still the case (assuming that X is compact) that $\hat{\phi}$ is either identically zero or A-covariantly constant; and it is still the case that $a$ can be written as $\nu_{\ddagger}\hat{\phi}$ with $\nu_{\ddagger}$ being a self-dual, harmonic 2-form. However, $F_A^+$ is no longer zero; it is $m\nu_{\ddagger}\hat{\phi}$.

The preceding conclusions also hold when X is non-compact and certain growth conditions on $|\hat{\phi}|$ are imposed. For example, if X is complete with subquadratic volume growth as a function of distance, and if $|\hat{\phi}|$ is bounded, then the conclusions of the preceding paragraph hold. (Subquadratic volume growth means that the volume as a function of distance r from any given point is $o(r^2)$ for r >> 1. Manifolds with cylindrical ends are examples because their volume grows as $\mathcal{O}(r)$.)

In any event, there is an analog of Theorems 1.1 and 1.2 for the equations in (1.6). These analogs are stated in Section 7b. Section 7b also points out how their proofs differ from those of Theorems 1.1 and 1.2. (The differences are all cosmetic.)

**a) Outline of the proof of Theorem 1.2**

If $(A, a)$ obeys (1.2), then $(A, a)$ obeys $d_A^\dagger d_A a = 0$. A Bochner Weitzenboch formula that rewrites the equation $d_A^\dagger d_A a = 0$ in terms of the covariant Laplacian (on $C^\infty(X; \Lambda^+ \otimes \mathcal{G})$ and an action of $F_A^+$ and an action of the Riemann curvature tensor on $a$. An integral version of this Bochner-Weitzenboch formula leads directly to an a priori bound for the integral of the function $|\nabla_A a|^2 + r^2|[a; a]|^2$ on open sets in X with compact closure. This leads in turn to an a priori $L^2$ bound for the norm of $d|a|$ because the connection A and the Levi-Civita connection are metric compatible. (Keep in mind that $|d|a|| \leq |\nabla_A a|$.) The Bochner-Weitzenboch formula also leads to an a priori pointwise bound for $|a|$. The Bochner-Weitzenboch formula and these simple consequences are discussed in Section 2. (Theorem 1.1's proof uses these consequence of the Bochner-Weitzenboch formula, Karen Uhlenbeck's theorems in [U1] and [U2], and some standard elliptic regularity theorems. The proof is sketched in Section 2b.)

With regards to Theorem 1.2: The observations in the preceding paragraph are applied to each term in the sequence $\{(A_n, a_n)\}_{n \in \mathbb{N}}$ and, with the third bullet in (1.2), they imply that the sequence of functions $\{|a_n|\}_{n \in \mathbb{N}}$ is bounded in the $L^2_1$ topology on open sets in X with compact closure. As consequence, there is a subsequence of $\{|a_n|\}_{n \in \mathbb{N}}$ that converges on compact subsets of X in the weak $L^2_1$ topology. The resulting limit function is the function $|\nu_{\ddagger}|$. The aforementioned Bochner-Weitzenboch formula also leads (using the maximum principle) to an a priori pointwise bound for each $|a_n|$. The existence of these bounds imply that $|\nu_{\ddagger}|$ is an $L^\infty$ function. Proposition 6.1 makes a formal assertion of what was just said with regards to $L^2_1$ convergence and the $L^\infty$ bound. It also gives a useful pointwise definition for $|v_\ddagger|$. This same proposition also specifies Theorem 1.2's subsequence $\Lambda$.

The four bullets of Theorem 1.2 are proved with the help of certain a priori estimates and bounds for the $\Lambda^+ \otimes \mathfrak{G}$ component of any given solution to (1.2). Let (A, $a$) again denote a solution. These bounds are, for the most part, local bounds on balls in X for the integral of $|a|^2$ on the boundary of the ball and for certain integrals of the function $|\nabla_A a|^2 + r^2 |[a;a]|^2$ over the interior. In particular, the square root of the averages of $2\pi^2 |a|^2$ on the boundaries of balls centered at a given point are used to define a positive function on an interval of the form $(0, r_0)$ with $r_0$ being a fixed number that is determined by the local geometry of X. This function is denoted by $\kappa$. (The precise definition modifies the square root by a multiplicative factor that differs from 1 by $\mathcal{O}(r^2)$.) A second function, denoted by $\mathcal{N}$, is defined by writing the derivative of $\kappa$ at any given radius r as $\frac{d\kappa}{dr} = \frac{\mathcal{N}}{r} \kappa$. (As it turns out, the function $\mathcal{N}$ is the product of $\frac{1}{r^2 \kappa^2}$ with the integral of the function $|\nabla_A a|^2 + r^2 |[a;a]|^2$ over the ball of radius r.) Section 3 defines $\kappa$ and $\mathcal{N}$ and derives some consequences of the definitions. (The functions $\kappa$ and $\mathcal{N}$ are inspired by a pair of analogous functions that were introduced by Almgren [Al] many years ago to study singularities in minimal surfaces. These analogs have subsequently been used to great effect in differential geometry and analysis. See for example [DF], [HHL] and, very recently, [DeL].)

The behavior of $\kappa$ and $\mathcal{N}$ and what they say about (A, $a$) depends on whether the radius r is below or above a certain critical radius that is associated to the center point of the given ball. This critical radius is denoted by $r_{c_\mu}$. The number $r_{c_\mu}$ is defined (see Section 4) after choosing a number $\mu$ from $(0, \frac{1}{100}]$ and then a number $c$ which is large, with an $\mathcal{O}(1)$ lower bound that is determined a priori by only $\mu$ (and the geometry of X around the given point if X is not compact). The number $r_{c_\mu}$ for a given point p where $a \neq 0$ is the largest radius r of a ball centered at p such that the following condition holds: If q is any point in the concentric ball of radius $(1-\mu)r$ centered at p, then the integral over the radius r ball centered at p of the function

$$\frac{1}{\text{dist}(q, \cdot)^2} (|\nabla_A a|^2 + r^2 |[a;a]|^2)$$

(1.7)

is no greater than $\frac{1}{c^2}$. If $a = 0$ at p, then $r_{c_\mu}$ is set equal to 0. (If $a = 0$ at the point, then there is no r where the q = p version of (1.7) is less than $\frac{1}{c^2}$ if c is greater than $\mathcal{O}(1)$.)

The number $r_{c_\mu}$ is significant in two ways: First, if r is less than $r_{c_\mu}$, then the point-wise norm of $|a|$ in the radius $(1-\mu)r$ concentric ball differs from $\frac{1}{\sqrt{2\pi}} \kappa(r)$ by very

much less than $\kappa(r)$. This implies that the pointwise norm of $|a|$ is nearly constant in this slightly smaller ball. In addition, if r is such that $rr\kappa(r)$ is not too small (too small being $\mathcal{O}(1)$), then $a$ can be written on the radius $(1-\mu)r$ ball as $a = v\sigma_\ddagger + \mathfrak{a}$ with $v$ being a self-dual 2-form that is nearly closed, with $\sigma_\ddagger$ being a unit length section of $\mathcal{G}$ that is nearly A-covariantly constant, and with $\mathfrak{a}$ being a section of $\Lambda^+ \otimes \mathcal{G}$ that is very small pointwise and has small A-covariant derivative. Small here means $\mathcal{O}(e^{-x})$ with $x$ being $\mathcal{O}(\sqrt{rr\kappa})$. This decomposition when applied to each $n \in \Lambda$ version of $(A_n, a_n)$ leads to corresponding $v_n$ and $\sigma_n$ with $\sigma_n$ as described by Theorem 1.2 and with $v_\ddagger$ from Theorem 1.2 given by a limit of the sequence $\{v_n\}_{n \in \Lambda}$.

What happens if $r \geq r_{c\mu}$? What happens here is the following (and this is really the key to the whole proof): If $r > r_{c\mu}$, then the value of $N(r)$ is bounded uniformly away from zero. This is to say that there exists a positive number determined only by $c$ and $\mu$ (denoted here by $\delta$) such that $N(r) > \delta$ if $r > r_{c\mu}$. This implies in turn (because of $N$'s relation to the derivative of $\kappa$) that $\kappa(r) \leq (\frac{r}{s})^\delta \kappa(s)$ if $r < s$ and both are greater than $r_{c\mu}$ (and less than a number that is $\mathcal{O}(1)$). Because $\kappa$ is an increasing function (its derivative is positive), the latter bound implies that $r_{c\mu}$ can not be too small if $|a|$ at the center of the ball is not too small. Thus, one has an apriori lower bound for $r_{c\mu}$. This lower bound and what is said in the preceding paragraph about when $r < r_{c\mu}$ leads in a straightforward way to all of the conclusions of Theorem 1.2. (How this happens is explained in Section 6.)

The proof that $N(r) > \delta$ if $r > r_{c\mu}$ is surprisingly straightforward. In particular, the proof requires nothing by way of a priori bounds for the anti-self dual part of the curvature of A. (More is said momentarily about what the proof doesn't use.) What the proof does use (and this is crucial) is the fact that the equation $d_A a = 0$ in (1.2) is a first order, elliptic equation; and in particular, the fact that the A-covariant derivative of any given $\Lambda^+$ component of $a$ (when written using a local frame for $\Lambda^+$) is determined pointwise by the other components of $a$. The proof also uses the fact that the term $r^2[a;a]$ that appears in (1.2) is involves commutators of $a$.

The definition of $r_{c\mu}$ and the analysis of $N$ where $r$ is greater than $r_{c\mu}$ is in Section 4 of this paper. The analysis of $a$ where $r < r_{c\mu}$ is in Section 5.

A formally similar strategy was used in [T2] and [T3] (and [T4] and also [HW], both of which concern 3-dimensional manifolds) to prove analogs of Theorem 1.2 for equations that have the same schematic form as the equations in (1.2). However, the implementation of the strategy is very different in these papers and, as a consequence, so is much of their technology. There are two salient differences: First, the analog of $r_{c\mu}$ is defined in these papers in terms of a threshold size for the $L^2$ norm of the full curvature $F_A$ on balls centered at the given point. Second, the analysis to characterise the behavior of K where r is greater than the $r_{c\mu}$ analog required a detailed understanding of the r-derivative of the function N. In particular, the basic method of attack was to prove a 'monotonicity' formula for N to the effect that once N is larger than some small positive number, then it never gets smaller as r increases. (The monotonicity formula in these references was inspired by a monotonicity formula for an analog of N that was introduced by Almgren [Al]. A similar monotonicity formula is exploited in all of the other analogs of N that have appeared over the years.) A great deal of work was done in [T2], [T3] (and [T4] and [HW]) to prove this monotonicity of N directly from the formula for its derivative (this was the method of proof used by Almgren for the ancestral incarnation of the function N.) These monotonicity proofs required (or seem to require) a priori bounds for the integrals of $|F_A|$ and $|F_A|^2$ on balls of radius greater than the $r_{c\mu}$ analog. The author doesn't know any analogous local curvature bounds for solutions to (1.2). (There are global bound for the whole of X when X is compact that come from the Pontrjagin class integral.)

In any event, it turns out that neither a priori bounds for the whole of $F_A$ nor the complicated formulas in [T2] and [T3] for the derivative of N are needed to prove that N is greater than a certain canonical positive number where $r > r_\delta$. (It is most likely the case that the theorems in [T2], [T2] and [T4] can be reproved with much less work using the ideas and associated technology that is introduced in this paper.)

**b) Conventions**

This subsection lays out various conventions (notational and otherwise) that are used in this article for the most part with no further comment.

ORTHONORMAL FRAMES: The derivation of the various algebraic identities (and their applications) is facilitated by working locally and writing $a$ using and its covariant derivatives using chosen orthonormal frame for T*X and $\Lambda^+$. To elaborate, let $\{e^\alpha\}_{\alpha=1,2,3,4}$ denote an oriented, orthonormal frame for T*X defined near some given point in X. (Greek indices run from 1 to 4.) A 1-form $s$ is written using $\{e^\alpha\}_{\alpha=1,2,3,4}$ as $s_\alpha e^\alpha$. By the same token, the four components of the Levi-Civita connection's covariant derivative on tensors are written as $\{\nabla_\alpha\}_{\alpha=1,2,3,4}$; and the four components of the covariant derivative on $\mathfrak{G}$-valued tensors defined by a connection A on P and the Levi-Civita connection are written as $\{\nabla_{A\alpha}\}_{\alpha=1,2,3,4}$. A corresponding orthonormal frame for $\Lambda^+$ is always used:

$$\omega^1 = \tfrac{1}{\sqrt{2}}(e^2 \wedge e^3 + e^1 \wedge e^4), \quad \omega^2 = \tfrac{1}{\sqrt{2}}(e^3 \wedge e^1 + e^2 \wedge e^4), \quad \omega^3 = \tfrac{1}{\sqrt{2}}(e^1 \wedge e^2 + e^3 \wedge e^4).$$
(1.8)

Supposing that $t$ denotes a section of $\Lambda^+$ (or the tensor of $\Lambda^+$ with $\mathfrak{G}$ or some other bundle), then it can be written using the frame in (1.8) as $t = t_a\omega^a$. These Latin subscripts run from 1 to 3.

INNER PRODUCTS: The Lie algebra of SU(2) (which is the fiber of the bundle $\mathfrak{G}$) is identified henceforth with the vector space of $2 \times 2$, traceless, anti-Hermitian complex matrices. This vector space is denoted henceforth at $\mathfrak{s}$. The fiber norm on $\mathfrak{G}$ was defined using an ad-invariant inner product on $\mathfrak{s}$ which is taken to be the following one: Supposing that $\mathfrak{a}, \mathfrak{b}$ are elements in $\mathfrak{s}$, then their inner product is $-\frac{1}{2}\text{trace}(\mathfrak{ab})$. By way of notation, if $\mathfrak{a}$ is any $2 \times 2$ complex matrix, then $-\frac{1}{2}\text{trace}(\mathfrak{a})$ is denoted by $\langle\mathfrak{a}\rangle$. Thus, the inner product of $\mathfrak{a}$ and $\mathfrak{b}$ is $\langle\mathfrak{ab}\rangle$. With regards to inner products in general, if $\mathfrak{a}$ and $\mathfrak{b}$ are tensors or $\mathfrak{G}$-valued tensors, then their inner product (which is denoted by $\langle\mathfrak{a},\mathfrak{b}\rangle$) is defined using the preceding inner product for $\mathfrak{G}$ and the Riemannian inner product on tensors. The pointwise norm of a section is defined using this inner product. If $\mathfrak{a}$ is the section, then the norm is denoted by $|\mathfrak{a}|$. For example, if $\mathfrak{a}$ is sections of $\Lambda^+ \otimes \mathfrak{G}$, then its norm (when the section is written using a basis as in (1.8)) is $|\mathfrak{a}| = \langle\mathfrak{a}_a\mathfrak{a}_a\rangle^{1/2}$.

WHEN X IS NON-COMPACT: Theorem 1.2 is really a local theorem about the behavior of solutions to (1.2) on open subsets of X. In particular, no assumptions are stated about the global behavior of the metric. However, to facilitate various applications of integration by parts, it is necessary (when X is not compact) to fix a priori an open subset in X with compact closure and to work implicitly only on this subset. This makes it possible to state theorems, propositions, lemmas and inequalities with constants that are position independent. The assumption that all activity takes place on this chosen subset is implicit (for the most part) in the statement and proofs of all of the assertions that follow. The chosen subset is denoted by $\underline{X}$ in what follows.

CONSTANTS INDEPENDENT OF $r$ and $(A, a)$: The convention in this paper is to use $c_0$ to denote a number that is greater than 1 and independent of the value of $r$ in (1.2) and the any given solution $(A, a)$ to (1.2). The precise value of $c_0$ can be assumed to increase between successive incarnations. Given a number $\mu$ (which is chosen from $(0, \frac{1}{100}]$), the convention has $c_\mu$ denoting a number that is greater than 1 and indpendent of $r$ and $(A, a)$, but it can depend on the choice of $\mu$.

CUT-OFF FUNCTIONS: Cut off functions are used that are equal to 1 on a given ball in X and equal to 0 on the complement of a larger, concentric ball. These can all be built from a fiducial function on $\mathbb{R}$ to be denoted by $\chi$ that is non-increasing and equal to 1 on $(-\infty, \frac{1}{2}]$ and equal to 0 on $[\frac{3}{4}, \infty)$. By way of an example, if the center point of the balls is p and the radius of the inner ball is r and that of the outer ball is r+h, then the function $\chi(\frac{1}{h}(\text{dist}(\cdot,p) - r))$ has the desired properties (here, dist($\cdot$,p) is the function on X that gives the distance to p.) Note that the norm of the derivative of this function is bounded by $c_0 \frac{1}{h}$, and the norm of its covariant Hessian is bounded by $c_0 \frac{1}{h^2}$. (Remember the convention about $c_0$!)

**c) Basic function space inequalities and Green's functions**

As was the case in [T2] and [T3], certain 4-dimensional function space inequalities will be invoked. Also, Green's functions are used. More is said momentarily. To set the notation, fix p ∈ X and given r ∈ (0, $c_0^{-1}$] with the upper bound here being far less than the radius of a Gaussian coordinate chart centered at p. Use $B_r$ to denote the open, radius r ball centered at p. The boundary of the closure of $B_r$ is denoted by $\partial B_r$; it is the sphere of radius r centered at p.

FUNCTION SPACE INEQUALITIES: Let $L^2_1$ denoting the space of square integrable functions on $B_r$ with square integrable differential. If $f$ is an $L^2_1$ function on $B_r$, then $f^2$ and $\text{dist}(\cdot,p)^{-1} f$ are square integrable on $B_r$; and their norms obey the following inequalities:

- $(\int_{B_r} f^4)^{1/2} < c_0 \int_{B_r} (|df|^2 + \frac{1}{r^2}|f|^2)$.

- $\int_{B_r} \frac{1}{\text{dist}(p,\cdot)^2} f^2 < c_0 \int_{B_r} (|df|^2 + \frac{1}{r^2}|f|^2)$.

- $\int_{B_r} f^2 \leq c_0 (r^2 \int_{B_r} |df|^2 + r \int_{\partial B_r} |f|^2)$

(1.9)

The first inequality is a Sobolev inequality and the second is called Hardy's inequality. The third is proved using an integration by parts. (See, e.g. [AF] for a general reference on Sobolev inequalities.) Note that if $s$ is a tensor on $B_r$, then the norm of d|s| is no greater than |∇s|; and if $s$ is a 𝔊-valued tensor, then the norm of d|s| is not greater than $|\nabla_A s|$ for any connection A on P. This implies in particular, that the inequalities in (1.9) hold with $f$ replaced by |s| and with d$f$ replaced by ∇s or $\nabla_A s$ as the case may be.

GREEN'S FUNCTIONS: Green's functions for the Dirichlet Laplacian on balls in X play role also in the proof (as they do in [T2] and [T3]). To say more about these functions, suppose that p is a given point in X and that r > 0 is much less than the injectivity radius at p. The Laplacian on $B_r$ is the operator $d^\dagger d$ acting on smooth functions. Here, $d^\dagger$ is the formal $L^2$ adjoint of d. Given a point q ∈ $B_r$, there is a Green's function for the Laplacian which is smooth on $B_r$–q, positive on $B_r$ and zero on $\partial B_r$. This Green's function is denoted by $G_q$ and it is described by the following:

- $G_q \leq c_0 \frac{1}{\text{dist}(q,\cdot)^2}$ .
- $|\nabla G_q| \leq c_0 \frac{1}{\text{dist}(q,\cdot)^3}$ .
- $G_q \geq c_0^{-1} \frac{1}{\text{dist}(q,\cdot)^2}$ *where the distance to* q *is less than* 4 *times the distance to* $\partial B_r$.

(1.10)

By way of a parenthetical remark (for now): If $f$ is an $L^2_1$ function on $B_r$, then it is a consequence of the second bullet in (1.9) and the first bullet in (1.10) that the function $f^2 G_p$ is integrable on $B_r$.

### 2) Bochner-Weitzenboch formula

A Bochner-Weitzenbock formula for the equations in (1.2) plays a central role in Mare's proof of Theorem 1.1 and in the upcoming proof of Theorem 1.2. Section 2a exhibits this Bochner-Weitzenbock formula and derives some immediate consequences (Lemma 2.1). Lemma 2.1 is then used in the sketch given in Section 2b of Ben Mare's proof of Theorem 1.1.

### a) The Bochner-Weitzenboch formula

Let $*$ denote the metric's Hodge star operator. This maps 3-forms to 1-forms. The equation in (1.2a) saying $d_A a = 0$ implies in turn that $d_A * d_A a = 0$ which is a $\mathfrak{G}$-valued 2-form identity. The self-dual part of this equation can be written (after commuting covariant derivatives) and using (1.2b) as

$$\nabla_A^\dagger \nabla_A a + r^2 [a_c, [a, a_c]] + \mathfrak{R} a = 0$$

(2.1)

where the notation has $\nabla_A^\dagger$ denoting the formal, $L^2$ adjoint of the covariant derivative $\nabla_A$; and it has $\mathfrak{R}$ denoting an endomorphism of $\Lambda^+$ (promoted to one of $\Lambda^+ \otimes \mathfrak{G}$) that is linear in the metric's scalar curvature and the self-dual part of the metric's Weyl curvature. Of particular note here is that $\mathfrak{R}$ and its $\nabla$-covariant derivatives to any given order are uniformly bounded.

Taking the inner product of (2.1) with $a$ leads to the following equation:

$$\tfrac{1}{2} d^\dagger d |a|^2 + |\nabla_A a|^2 + r^2 |[a;a]|^2 + \langle a, \mathfrak{R} a \rangle = 0.$$

(2.2)

(The $|[a;a]|^2$ term is the inner product of $a$ with $[a_c,[a,a_c]]$. The proof that this is so uses the cyclic property of the trace $\langle \cdot \rangle$ and the fact that $|[a;a]|^2 = \langle [a_a,a_b][a_a,a_b]\rangle$.)

The identity in (2.2) leads in turn to an analog of [T3]'s Lemma 2.1:

**Lemma 2.1**: *Fix an open subset* $\underline{X} \subset X$ *with compact closure and there exists* $\kappa > 1$ *with the following significance: Fix* $r > 0$ *and let* $(A, a)$ *denote a solution to (1.2) on X.*

- $\int_{\underline{X}} (|\nabla_A a|^2 + r^2 |[a;a]|^2) \leq \kappa$.

- $\sup_{p \in \underline{X}} |a|(p) \leq \kappa$.

- *If* $p \in \underline{X}$, *then* $\int_{\underline{X}} \frac{1}{\text{dist}(p,\cdot)^2} (|\nabla_A a|^2 + r^2 |[a;a]|^2) \leq \kappa$.

*Proof of Lemma 2.1*: The proof is identical but for notation to the proof of Lemma 2.1 in [T2] when X is compact and almost the same when X isn't: One first multiples (2.2) by a suitable nonnegative function (with compact support when X is not compact) and integrates over X. Integration by parts is then used to rewrite the resulting integral identity. Use a compactly supported function that is equal to 1 on $\underline{X}$ for the first bullet. For the second and third bullets, use a function with compact support on a small radius ball in X centered at p that is equal very near p to the Dirichelet Green's function for the metric Laplacian $d^\dagger d$ with pole at p

**b) Proof of Theorem 1.1**

This subsection sketches Ben Mare's proof of Theorem 1.1 (which is in [Ma]). By way of terminology: An $L^2_{1;loc}$ connection on P over a given open set in X can be written as $A_0 + \mathbb{A}$ where $A_0$ is a smooth connection defined on the open set, and where $\mathbb{A}$ is a section of $T^*X \otimes \mathfrak{G}$ on the open set whose components with respect any product structure over any ball with compact closure in the open set are from the Sobolev space $L^2_1$ (the norms of the components and their first derivatives are square integrable on the ball in question).

To start the proof, suppose for the moment that $(A, a)$ obeys (1.2) and that X is compact. By virtue of the top bullet of Lemma 2.1, the square of the $L^2$ norm of $F_A^+$ is bounded by $c_0 r^2$. This implies in turn that the square of the $L^2$ norm of $F_A^-$ is bounded by $c_0(1 + r^2)$ because the difference between these numbers is proportional to the first Pontrjagin number of the bundle P.

With the preceding understood, suppose that $\{r_n, (A_n, a_n)\}_{n \in \mathbb{N}}$ is a sequence of solutions to (1.2) with $\{r_n\}_{n \in \mathbb{N}}$ being bounded. If X is not compact, assume in addition that the sequence whose n'th component is the $L^2$ norm on X of the curvature of $A_n$ is also bounded. Given the preceding, then Uhlenbeck's theorem [U1] can invoked to find the following data:

- *A finite set of points set of points in* X. *The set is denoted by* Ξ.
- *An* $L^2_{1;loc}$ *connection on* P *over* X−Ξ *to be denoted by* $A_\ddagger$.
- *A subsequence* $\Omega \subset \mathbb{N}$
- *A sequence of automorphisms of* P *indexed by* Ω. *These are denoted by* $\{g_n\}_{n\in\Omega}$.

(2.3)

This data is such that the sequence $\{g_n{}^*A_n\}_{n\in\Omega}$ converges to $A_\ddagger$ in the weak $L^2_1$ topology on any open set in X−Ξ with compact closure.

It follows from the preceding observation and the $r = r_n$ and $(A, a) = (A_n, a_n)$ version of Lemma 2.1 that the sequence $\{g_n{}^*a_n\}_{n\in\Omega}$ is bounded in the $L^2_1$ topology on any open set in X−Ξ with compact closure. As a consequence, there is an $L^2_{1;loc}$ section of $\Lambda^+\otimes\mathfrak{G}$ on X−Ξ (it is denoted by $a_\ddagger$) and a non-negative number $r'$ and a subsequence $\Lambda \subset \Omega$ such that $\{r_n\}_{n\in\Lambda}$ converges to $r'$ and such that $\{g_n{}^*a_n\}_{n\in\Lambda}$ converges in the weak $L^2_{1;loc}$ topology to $a_\ddagger$

The preceding implies that $(A_\ddagger, a_\ddagger)$ obeys the top two bullets of (1.2) with $r = r'$. The third bullet of (1.2) is also obeyed by $a_\ddagger$ if X is compact. Otherwise, it can be said that the $L^2$ norm on X of $a_\ddagger$ not greater than 1. Because the top two bullets of (1.2) are obeyed, there is necessarily an automorphism of the bundle P over X−Ξ that makes $(A_\ddagger, a_\ddagger)$ smooth. (This is a standard elliptic regularity argument.)

Granted that $(A_\ddagger, a_\ddagger)$ is smooth, the convergence described above can be revisited to see that the automorphisms $\{g_n\}_{n\in\Omega}$ from (2.3) can be chosen so that the sequence $\{(g_n{}^*A_n, g_n{}^*a_n)\}_{n\in\Theta}$ converges to $(A_\ddagger, a_\ddagger)$ in the $C^\infty$ topology on compact subsets of X−Ξ.

The final step in the proof revisits and slightly modifies Uhlenbeck's arguments in [U2] to prove that there is a principal bundle $P' \to X$ and an isomorphism from $P'|_{X-\Xi}$ to $P|_{X-\Xi}$ that identifies $(A_\ddagger, a_\ddagger)$ with a smooth pair of connection on P' that is defined over all of X (this is A´) and section of $\Lambda^+\otimes\mathfrak{G}'$ that is also defined on the whole of X (this is $a'$). Since these are smooth, they obey the top two bullets of (1.2) on the whole of X.

## 3. The functions κ, N

Fix a point $p \in X$ and fix a positive number $r_0$ so that the ball of radius $r_0$ centered at p is well inside a Gaussian coordinate chart centered at p. (If X is compact or if $\underline{X}$ has been chosen, then $r_0$ can be taken to be $\mathcal{O}(c_0^{-1})$.) Two functions on $(0, r_0]$ are going to be used to describe the behavior of the norm of $a$ on the ball of radius $r_0$ centered at p.

These are K and N. Version of these functions are also introduced for each component of
*a* as defined by a suitably chosen orthonormal frame for $\Lambda^+$ on the ball $B_{r_0}$.

**a) The definition of K and N**

The definitions that follow for K and N and the subsequent discussion differs only in notation from what is said in Section 3b of [T3] (and the latter borrows from Section 3a in [T2].) Because of this, the discussion below is a summary only. The reader is directed to Section 3b of [T3] (see also Section 3a of [T2]) for proofs of what is said.

The definition of K requires the introduction of two other functions on $(0, r_0]$ (which are used only to get to K). The first is the function h that is defined by the rule

$$r \to h(r) = \int_{\partial B_r} |a|^2 \ .$$

(3.1)

An appeal can be made to theorem of Aronszajn [Ar] to prove that $h(r) > 0$ for $r > 0$ unless *a* is identically zero. (One can also prove this by borrowing the argument in Section 2c of [T5].) The second function is called $\mathfrak{d}$ and it is defined by the rule:

$$r \to \mathfrak{d}(r) = \int_0^r \left( \frac{1}{h(s)} \left( \int_{B_s} \langle a, \mathfrak{R}a \rangle + \tfrac{1}{2} \int_{\partial B_s} M|a|^2 \right) \right) ds \ .$$

(3.2)

where $\mathfrak{R}$ is the curvature endomorphism in (2.1) and where M is defined by writing the second fundamental form of $\partial B_r$ for any given $r \in (0, r_0]$ as $\tfrac{3}{r} + M$. Lemma 2.1 implies that $|\mathfrak{d}(r)| \le c_0 r^2$.

With h and $\mathfrak{d}$ in hand, define the desired function K by the rule

$$K = e^{\mathfrak{d}} \left( \tfrac{1}{r^3} h \right)^{1/2}$$

(3.3)

Up to an $\mathcal{O}(r^2)$ factor, this function K is defined so that its square is $2\pi^2$ times the average of $|a|^2$ on the radius r sphere centered at p. Meanwhile, the desired function N is defined by the rule

$$r \to N(r) = \frac{1}{r^2 K(r)^2} \int_{B_r} (|\nabla_A a|^2 + r^2 \,|[a; a]|^2) \ .$$

(3.4)

The value of N at any given $r \in (0, r_0]$ can also be written as

$$N(r) = \frac{1}{r^2 K(r)^2} \int_{B_r} (|\nabla_A a|^2 + 2r^{-2} |F_A^+|^2)$$

(3.5)

because of the second bullet in (1.2). As was the case in [T2] and [T3], the function ∂ is introduced in the definition of K so that N and K are related via the identity

$$\frac{d}{dr} K = \frac{1}{r} NK .$$

(3.6)

Thus, if r and s are from $(0, r_0]$ and $r > s$, then $K(r)$ and $K(s)$ are related by the rule

$$K(r) = \exp\left(\int_s^r \frac{N(\tau)}{\tau} d\tau\right) K(s)$$

(3.7)

(The formula in (3.6) is proved by first taking the inner product of both sides of (2.1) with a. Having done so, integrate the result over $B_r$. The boundary term differs from $\frac{1}{2} r^3 \frac{d}{dr} K^2$ by an $\mathcal{O}(r^2 K^2)$ factor that is accounted for by one of the terms in the derivative of the function ∂. Meanwhile, the integral over $B_r$ that results from the integration by parts can be rewritten using the definition of ∂ so as to obtain (3.6) with the preceding identification of the boundary term.)

Note that (3.6) and (3.7) imply that K is an increasing function on $(0, r_0]$. This implies in particular that if $r \in (0, r_0]$ and $s \in [0, r]$, then

$$(1-c_0 r^2) \frac{\pi^2}{2} (r^4 - s^4) K^2(s) \leq \int_{B_r - B_s} |a|^2 \leq (1+c_0 r^2) \frac{\pi^2}{2} (r^4 - s^4) K^2(r) .$$

(3.8)

Although K can not be zero at positive r, it is none-the-less possible for N to get very large. Even so, one has the following: Supposing that $\mu \in (0, 1)$, then

$$(NK^2)|_{(1-\mu)r} \leq c_\mu K^2(r) .$$

(3.9)

This bound is derived by first multiplying both sides of (3.6) by K to get a formula for the derivative of $K^2$ saying that $\frac{d}{dr} K^2 = \frac{2}{r} NK^2$. Then, integrate the latter from $(1-\mu)r$ to r. The bound in (3.9) follows directly from this integral identity because the function $r \to r^2 K^2(r) N(r)$ is positive and increasing (it is the integral on $B_r$ of $|\nabla_A a|^2 + r^2 |[a; a]|^2$).

b) **The pointwise norm of a**

The values of the functions K and N on balls of radius $r \in (0, r_0]$ centered a p can used to bound the pointwise norm of |a| on slightly smaller radius, concentric balls. This

is what is said in effect by the following proposition (it the Vafa-Witten analog of Proposition 3.1 in [T3]).

**Proposition 3.1**: *There exists $\kappa > 1$; and given $\mu \in (0, \frac{1}{100}]$, there exists $\kappa_\mu > \kappa$; these numbers having the following significance: Fix $r > 1$ and suppose that $(A, a)$ is a solution to (1.2). Fix $p \in X$ (in $\underline{X}$ if X is non-compact) to define the functions K and N. Supposing that $r \in (0, \kappa^{-1}]$, then*

- $|a| \le \kappa_\mu K(r)$ *on* $B_{(1-\mu)r}$.

*Moreover, if $N(r) \le 1$, then*

- $|a| \le (1 + \kappa_\mu(r + \sqrt{N(r)})) \frac{1}{\sqrt{2\pi}} K(r)$ *on* $B_{(1-\mu)r}$.

*Proof of Proposition 3.1*: Except for notation, the proof is the same as the proof in Section 3c of [T3] for the Proposition 3.1 in [T3]. Even so, the proof is sketched because some of the arguments and steps are used again.

What follows gives a sketch of the proof of the top bullet of Proposition 3.1: Introduce the function $\chi$ from Section 1b and use it to define a function on $B_r$ that is equal to 0 where $\text{dist}(\cdot, p) \ge (1 - \frac{5}{8}\mu)r$, equal to 1 where $\text{dist}(\cdot, p) \le (1 - \frac{3}{4}\mu)r$. Denote this function by $\chi_\mu$. It can and should be constructed so that its derivatives obey the bounds $|d\chi_\mu| \le c_0 \frac{1}{\mu r}$, and $|\nabla d\chi_\mu| \le c_0 \frac{1}{\mu^2 r^2}$. Fix a point $q \in B_{(1-\mu)r}$ and introduce the Dirichelet Green's function with pole at q for $d^\dagger d$. This is the function $G_q$ from Section 1c. Multiply both sides of (2.1) by $\chi_\mu G_q$ and then integrate by parts. The result can be written schematically as

$$\tfrac{1}{2}|a|^2(q) + \int_{B_r} \chi_\mu G_q(|\nabla_A a|^2 + r^2|[a;a]|^2) = \int_{B_{(1-\mu/16)r}} L_{\mu,q}|a|^2 - \int_{B_{(1-\mu/16)r}} \chi_\mu G_q \langle a, \mathfrak{R}a\rangle$$

(3.10)

with $L_{\mu,q}$ denoting here and subsequently the function $-G_q d^\dagger d\chi_\mu + 2\langle d\chi_\mu, dG_q\rangle$.

The absolute value of the integral of $G_q\langle a, \mathfrak{R}a\rangle$ that appears on the right hand side of (3.10) is less than $c_0 r^2 K^2(r)$. This is proved by first bounding the integrand using the bound in the top bullet of (1.10) by $c_0 \frac{1}{\text{dist}(q,\cdot)^2}|a|^2$. Having done this, invoke the version of the middle bullet of (1.9) with p replaced by q and r replaced by $\mu r$. This leads to a bound on the $G_q\langle a, \mathfrak{R}a\rangle$ integral in (3.10) by $c_\mu$ times the integral over $B_{(1-\mu/64)r}$ of the function $|\nabla_A a|^2 + \frac{1}{r^2}|a|^2$. The next to last step invokes (3.9) (with $\mu$ replaced by $\frac{1}{64}\mu$) to bound the integral of $|\nabla_A a|^2$ over $B_{(1-\mu/64)r}$ by $c_\mu r^2 K^2(r)$; and the last step invokes the $s = 0$ version of the left most inequality in (3.8) to bound the integral over $B_r$ of $|a|^2$ by $c_\mu r^4 K^2(r)$ (and hence $\frac{1}{r^2}$ times this integral by $c_\mu r^2 K^2(r)$).

A bound for the absolute value of the $L_{\mu,q} |a|^2$ integral in (3.9) uses the first of the following two properties of the function $L_{\mu,q}$ when $q \in B_{(1-\mu)t}$. The second property is used momentarily in the proof of the second bullet of Proposition 3.1.

- $|L_{\mu,q}| \le c_\mu \frac{1}{r^4}$.
- $\int_{B_{(1-\mu/16)r}} L_{\mu,q} = 1$.

(3.11)

The top bullet follows from the top two bullets of (1.10); and the lower bullet follows via an integration by parts because $G_q$ is the Green's function for $d^\dagger d$ and the constant function 1 is annihilated by this same operator.

With the top bullet in (3.11) understood, then it follows that the $L_{\mu,q} |a|^2$ integral in (3.9) is bounded by $c_\mu \frac{1}{r^4}$ times then integral of $|a|^2$ over $B_r$, which is less than $c_\mu \kappa^2(r)$ because of what is said by the $s = 0$ version the of left most inequality in (3.8)).

The second bullet in (3.9) is proved by taking more care with the $L_{\mu,q} |a|^2$ integral in (3.10). To say more, view $\kappa(r)$ for the moment as a constant (so $d\kappa(r) = 0$). With this understood, replace $L_{\mu,q} |a|^2$ in (3.10) by $L_{\mu,q}(|a|^2 - \kappa^2(r))$ and then invoke the identity in the second bullet of (3.11) to obtain the inequality

$$\tfrac{1}{2}(|a|^2(q) - \tfrac{1}{2\pi^2}\kappa^2(r)) + \int_{B_r} \chi_\mu G_q (|\nabla_A a|^2 + r^2 |[a;a]|^2) \le \int_{B_{(1-\mu/4)r}} L_{\mu,q}(|a|^2 - \kappa^2(r)) + c_\mu r^2 \kappa^2(r).$$

(3.12)

(The $c_\mu r^2 \kappa^2(r)$ term on the far right accounts for the integral of $G_q \langle a, \Re a \rangle$ in (3.10).) The arguments in Part 2 from the proof of Proposition 3.1 in [T3] can be repeated almost word for word to bound the integral on the right hand side of (3.12) by $c_\mu \sqrt{N(r)} \kappa^2(r)$. Doing so leads from (3.12) to the bound

$$\tfrac{1}{2}(|a|^2(q) - \tfrac{1}{2\pi^2}\kappa^2(r)) + \int_{B_r} \chi_\mu G_q (|\nabla_A a|^2 + r^2 |[a;a]|^2) \le c_0(r^2 + \sqrt{N(r)}) \kappa^2(r).$$

(3.13)

This inequality with the inequality in Proposition 3.1's first bullet implies the inequality in Proposition 3.1's second bullet.

### c) The components of $a$

Fix $p \in X$ and then fix an oriented, orthonormal frame for $\Lambda^+|_p$. Parallel transport this frame along the geodesic arcs from p using the Levi-Civita connection to obtain an oriented, orthonormal frame for $\Lambda^+$ on the radius $r_0$ ball centered at p. Denote this frame by $\{\omega^c\}_{c=1,2,3}$ and write $a$ using this frame as $a = a_c\omega^c$. Each component is a section of the bundle $\mathfrak{G}$ over the radius $r_0$ ball centered at p.

It proves convenient to introduce a rough analog of the function K for $a_c$. This is the function on $(0, r_0]$ denoted by $K_c$ which is defined by the rule

$$K_c(r) = \left(\frac{1}{r^3} \int_{\partial B_r} |a_c|^2\right)^{1/2}.$$

(3.14)

(Lacking the $e^\partial$ factor, it is not as sophisticated as K.) There is a corresponding version of N that is denoted by $N_c$ which is defined at values of r where $K_c > 0$ by

$$N_c(r) = \frac{1}{r^2 K_c(r)^2} \int_{B_r} (|\nabla_A a_c|^2 + r^2 \,|[a; a_c]|^2)$$

(3.15)

These two functions are related through an analog of (3.6):

$$\frac{d}{dr} K_c = \frac{1}{r} N_c K_c + \mathfrak{e}_c r$$

(3.16)

with $\mathfrak{e}_c$ being a function on $(0, r_0]$ that obeys the bound

$$|\mathfrak{e}_c| \leq c_0 (1 + \sqrt{N_c})(1 + \sqrt{N})\, K(r).$$

(3.17)

The proof of (3.16) and (3.17) starts with (2.1): Take the inner product of both sides with the basis 2-form $\omega^c$. After some rewriting, the result is the identity

$$\nabla_A^\dagger \nabla_A(a_c) + r^2 [a_b, [a_c, a_b]] + \Gamma_{cb} \nabla_A(a_b) + \Upsilon_{cb} a_b = 0$$

(3.18)

with $\Gamma$ obeying $|\Gamma| \leq c_0 \text{dist}(\cdot, p)$ and with $\Upsilon$ obeying $|\Upsilon| \leq c_0$. Both are determined a priori by the Riemannian metric. Now, take the inner product of both sides of (3.18) with $a_c$ and integrate the result over $B_r$. An integration by parts then gives the integral identity

$$\int_{\partial B_r} \langle a_c \nabla_{Ar}(a_c) \rangle = \int_{B_r} (|\nabla_A a_c|^2 + r^2 \,|[a; a_c]|^2) + \mathfrak{r}$$

(3.19)

with $\nabla_{Ar}$ denoting the directional covariant derivative along the unit normal vector to $\partial B_r$, and with $\mathfrak{r}$ denoting a function of r that obeys the bound

$$|\mathfrak{r}| \leq c_0 \int_{B_r} (r\,|a_c||\nabla_A a| + |a_c||a|) \ .$$

(3.20)

The left hand side of (3.19) differs from $\frac{1}{2} r^3 \frac{d}{dr} K_c^2$ by a term with norm bounded by $c_0 r^5 K_c^2$. Meanwhile, the integral on the right hand side of (3.14) is $r^2 K_c^2 N_c$. Thus, (3.19) leads to (3.16) if $\mathfrak{r}$ obeys the bound $\mathfrak{r} \leq c_0 r^4 (1 + \sqrt{N_c})(1 + \sqrt{N}) K_c(r) K(r)$. And, just such a bound follows from the right hand side of (3.20). (To see this, first use Holder's inequality to obtain integrals over $B_r$ of $|a_c|^2$, of $|a|^2$, and of $|\nabla_A a|^2$. Then invoke the third bullet of (1.9) to bound the integrals of $|a_c|^2$ and $|a|^2$; and then invoke the respective definitions of $N$ and $N_c$.)

Whereas $N_c$ is not defined where $K_c(r) = 0$, the product $N_c K_c^2$ is well defined at all values of r. (It is probably never the case that $K_c = 0$ at positive r unless $a_c$ is identically zero on $B_r$. Certainly, such is the case if the metric on X is real analytic in Gaussian coordinates.) Anyway, it follows from the definitions that $N_c K_c^2 \leq c_0 N K^2$. Because of this, multiplying both sides of (3.16) by $K_c$ gives an equation for the derivative of $K_c^2$ that is valid at all positive values of r. The latter equation leads to a bound that is an analog of the bound in (3.9) which says the following: Given $\mu \in (0, 1)$, then

$$\frac{1}{r^2} \int_{B_r} (|\nabla_A a_c|^2 + r^2 \,\|[a; a_c]\|^2) \leq c_\mu (K_c^2(\tfrac{1}{1-\mu} r) + r^2 N(\tfrac{1}{1-\mu} r) K^2(\tfrac{1}{1-\mu} r)) \ .$$

(3.21)

There is also a pointwise bound $|a_c|$ that is analogous to the bound in the top bullet of Proposition 3.1. This bound is stated below.

**Proposition 3.2**: *There exists $\kappa > 1$; and given $\mu \in (0, \frac{1}{100}]$, there exists $\kappa_\mu > \kappa$; these numbers having the following significance: Fix $r > 1$ and suppose that $(A, a)$ is a solution to (1.2). Fix $p \in X$ (in $\underline{X}$ if X is non-compact) to define the functions $K$ and $N$. Fix $r \in (0, \kappa^{-1}]$. Supposing that $q \in B_{(1-\mu)r}$, then*

- $|a_c|(q) \leq \kappa_\mu (K_c(r) + r\, K(r))$ .
- $\int_{B_{(1-\mu)r}} G_q (|\nabla_A a_c|^2 + r^2 \,\|[a; a_c]\|^2) \leq \kappa_\mu (K_c^2(r) + r^2 K^2(r))$ .

*on the concentric ball $B_{(1-\mu)r}$.*

*Proof of Proposition 3.2*: The argument is much the same as the argument that is borrowed from Section 3c of [T3] to prove the top bullet of Proposition 3.1. What follows is a brief account: Take the inner product of both sides of (3.18) with $a_c$ to obtain a differential equation that has the schematic form

$$\tfrac{1}{2} d^\dagger d |a_c|^2 + |\nabla_A(a_c)|^2 + r^2 |[a, a_c]|^2 = \mathcal{M}$$

(3.22)

with $\mathcal{M}$ on $B_r$ obeying $|\mathcal{M}| \leq c_0 r \, (|a||\nabla_A a| + |a|^2)$

Fix a point $q \in B_{(1-\mu)r}$, multiply both sides of (3.18) by $\chi_\mu G_q$ and integrate the result over $B_r$. An integration by parts then leads to the identity

$$|a_c|^2(q) + \int_{B_{(1-\mu/16)r}} \chi_\mu G_q (|\nabla_A a_c|^2 + r^2 \, |[a; a_c]|^2) = \int_{B_{(1-\mu/16)r}} L_{\mu,q} |a_c|^2 + \int_{B_{(1-\mu/16)r}} \chi_\mu G_q \mathcal{M} \; .$$

(3.23)

The right most term on the right hand side of (3.23) is bounded by $c_\mu r^2 \kappa^2$. This is a consequence the top bullet of Proposition 3.1 and (3.13). The left most integral on the right hand side of (3.23) can be bounded by $c_\mu (\kappa_c^2 + r^2 \kappa^2)$ by appealing first to the top bulet in (3.11) to bound $L_{\mu,q}$, and then to the second two bullets in (1.9) and to (3.17) and (3.9) with their versions of $\mu$ replaced by $c_0^{-1} \mu$. With these bounds, then (3.23) says that $|a_c|^2(q) \leq c_\mu (\kappa_c^2 + r^2 \kappa^2)$ which implies what is said by the top bullet of Proposition 3.2. It also implies directly what is said by the lower bullet of Proposition 3.2.

**4. The number $r_{c\mu}$**

Fix a point $p \in X$ (in $\underline{X}$ if X is noncompact). This section introduces a number associated to p to be denoted by $r_{c\mu}$. As indicated by the notation, the definition requires the specification of a number $\mu \in (0, \tfrac{1}{100})$ and a number $c \in (1, \infty)$. (Eventually both $\mu$ and $c$ will be fixed to be $\mathcal{O}(1)$.) Having chosen $\mu$ and $c$ (and the point p), define $r_{c\mu}$ to be the largest number in $(0, r_0]$ obeying

$$\sup_{q \in B_{(1-\mu)r}} \frac{1}{K(r)^2} \int_{B_r} \frac{1}{\text{dist}(q,\cdot)^2} (|\nabla_A a|^2 + r^2 \, |[a,a]|^2) \leq \frac{1}{c^2} \; .$$

(4.1)

assuming that such a number exists. (Note that $\kappa$ in (4.1) is defined using p, the center point of $B_{(1-\mu)r}$.) If there is no number in $(0, r_0]$ where (4.1) holds, then set $r_{c\mu} = 0$. Since $\lim_{r \to 0} \frac{1}{\sqrt{2\pi}} K(r) = |a|(p)$ and since $|a|$ is continuous, the number $r_{c\mu}$ is strictly positive in the event that $a(p) \neq 0$. (The upcoming Proposition 5.1 implies the converse to this assertion: If $r_{c\mu} > 0$, then $a(p) \neq 0$.) Note in any event, the assignment of $r_{c\mu}$ to p varies continuously with variations of p.

The following proposition talks about the function $K$ at $r \geq r_{c_\mu}$. This proposition plays a very central role in the proof of Theorem 1.2. It implies (among other things) that $K(r)$ decreases as a positive power of $\frac{r}{r_0}$ as $r$ is decreased from $r_0$ as long as $r \geq r_{c_\mu}$.

**Proposition 4.1**: *There exists $\kappa > 1$, and given $\mu \in (0, \frac{1}{100}]$, there exists $\kappa_\mu > \kappa$, and given $c > \kappa_\mu$, there exists $\lambda_{c_\mu} > \kappa_\mu$; these numbers having the following significance: Fix $r > 1$ and suppose that $(A, a)$ is a solution to (1.2). Fix a point $p \in X$ (in $\underline{X}$ if $X$ is not compact) and use $\mu$, $c$ and $p$ to define $r_{c_\mu}$. Assume that $r_{c_\mu} < r < \lambda_{c_\mu}^{-1} r_0$ and that $r r K(r) \geq 1$. If $r_1 \in (r, \lambda_{c_\mu}^{-1} r_0]$ and if $r_2 \in [r, r_1]$, then $K(r_2) \leq (\frac{r_2}{r_1})^{1/\lambda_{c_\mu}} K(r_1)$.*

Proposition 4.1 is seen momentarily (in Section 4a) to follow almost as a corollary to the next proposition. This proposition makes an assertion that says (with minor qualifications) that $r_{c_\mu} \geq r$ if $N(r)$ is small.

**Proposition 4.2**: *There exists $\kappa > 1$, and given $\mu \in (0, \frac{1}{100}]$, there exists $\kappa_\mu > \kappa$, and given $c > \kappa_\mu$, there exists $\lambda_{c_\mu} > \kappa_\mu$; these numbers having the following significance: Fix $r > 1$ and suppose that $(A, a)$ is a solution to (1.2). Fix a point $p \in X$ (in $\underline{X}$ if $X$ is not compact) and use $\mu$, $c$ and $p$ to define the number $r_{c_\mu}$. Suppose that $r \leq \lambda_{c_\mu}^{-1} r_0$ but that $r r K(r) \geq 1$; and suppose in addition that $N(r) \leq \lambda_{c_\mu}^{-1}$. Then $r_{c_\mu} \geq r$.*

This proposition is proved in Section 4b.

a) **Proof of Proposition 4.1**

Take $\mu$ and $c$ to obey the constraints of Proposition 4.2; and take $\lambda_{c_\mu}$ from this proposition. Assume that $r < \lambda_{c_\mu}^{-1} r_0$ and that $r r K(r) \geq 1$. If $r_{c_\mu} < r$, then Proposition 4.2 asserts that $N \geq \lambda_{c_\mu}^{-1}$ on the whole of the interval $[r, \lambda_{c_\mu}^{-1} r_0]$. Granted this lower bound, then the claim made by Proposition 4.1 follows from (3.7).

b) **Proof of Proposition 4.2**

The proof has nine parts.

*Part 1*: The derivative of $N$ can be written as

$$\tfrac{d}{dr} N = \tfrac{1}{r^2 K(r)^2} \int_{\partial B_r} (|\nabla_A a|^2 + r^2 |[a; a]|^2) - \tfrac{2}{r} N(N+1).$$

(4.2)

Note in particular that if $N < \frac{1}{100}$, then $N$ can not decrease too fast as r increases because

$$\tfrac{d}{dr} N \geq -\tfrac{4}{r} N \quad \text{where} \quad N \leq \tfrac{1}{100}.$$

(4.3)

Therefore, if $r \in (0, r_0)$ and if $s \in (0, r]$, and if $N(r) \leq \frac{1}{100}$, then

$$N(r) \geq (\tfrac{s}{r})^4 N(s) \quad .$$

(4.4)

By way of a relevant example, this says in particular that if $\varepsilon \in (0, \frac{1}{100})$ and if $N(r) \leq \varepsilon$, then $N(s) \leq 16\varepsilon$ when $s \in [\tfrac{1}{2} r, r]$.

*Part 2*: Choose an oriented, orthonormal frame for $\Lambda^+$ at p. This frame is denoted by $\{\omega^1, \omega^2, \omega^3\}$. Then, fix an oriented, orthonormal frame for T*X at p so that (1.8) holds. (To do this, choose any unit length element in T*X at p to be $e^4$. Then set $e^c$ for $c \in \{1, 2, 3\}$ to be $*(e^4 \wedge \omega^c)$.) Parallel transport these frames along the radial geodesics from p using the Levi-Civita connection to obtain orthonormal frames for $\Lambda^+$ and T*X over the ball of radius $r_0$ centered at p.

For each $a \in \{1, 2, 3\}$ and each distinct pair $\alpha, \beta \in \{1, 2, 3, 4\}$, use $\eta^a{}_{\alpha\beta}$ to denote $\sqrt{2} *(\omega^a \wedge e^\alpha \wedge e^\beta)$. Thus, $\omega^a = \frac{1}{2\sqrt{2}} \eta^a{}_{\alpha\beta} e^\alpha \wedge e^\beta$ for each index a. The equation in (1.2) (which can be written as $\eta^a{}_{\alpha\beta} \nabla_\beta a_a = 0$) is equivalent to the equation

$$\nabla_{A\alpha} a_a + \varepsilon^{abc} \eta^c{}_{\alpha\beta} \nabla_{A\beta} a_c = 0$$

(4.5)

with $\{\varepsilon^{abc}\}_{a,b,c \in \{1,2,3\}}$ denoting the components of the completely anti-symmetric 3-tensor with $\varepsilon^{123} = 1$. This last identity leads to the observation that

$$|\nabla_A a_3|^2 \leq c_0 (|\nabla_A a_1|^2 + |\nabla_A a_2|^2) \quad .$$

(4.6)

(Of course, there are the two analogous identities that are obtained by switching index 1 with the indices 2 or 3.)

The inequality in (4.6) is really the key to the proof of Proposition 4.1 because of the following consequence: Supposing that $r \in (0, r_0]$, and supposing that $f$ is a non-negative, $L^1$ function on $B_r$, then

$$\int_{B_r} f(|\nabla_A a|^2 + r^2 \|[a; a]\|^2) \leq c_0 (\int_{B_r} f(|\nabla_A a_1|^2 + r^2 \|[a; a_1]\|^2) + \int_{B_r} f(|\nabla_A a_2|^2 + r^2 \|[a; a_2]\|^2)) .$$

(4.7)

Of particular interest for what follows are the cases when $f$ is the product of the Green's function $G_q$ for some given point $q \in B_{(1-\mu)r}$ and the characteristic function for $B_{(1-\mu)r}$. (This characteristic function is equal to 1 in $B_{(1-\mu)r}$ and it is equal to 0 on the complement of $B_{(1-\mu)r}$. It is denoted in what follows by $\theta_{r,\mu}$.)

*Part 3*: Given $q \in B_{(1-\mu)r}$, invoke (4.7) with $f = \theta_{t,\mu} G_q$. These various $q \in B_{(1-\mu)r}$ versions of (4.7) with the inequality in the second bullet of Proposition 3.2 imply that

$$\sup_{q \in B_{(1-\mu)r}} \int_{B_{(1-\mu)r}} G_q (|\nabla_A a|^2 + r^2 \, |[a;a]|^2) \leq c_\mu (K_1^2(r) + K_2^2(r) + r^2 K^2(r)) .$$

(4.8)

Therefore, if $r \leq c_\mu^{-1} \frac{1}{c}$ and if both $K_1(r)$ and $K_2(r)$ are less than $c_\mu^{-1} \frac{1}{c} K(r)$, then $r \leq r_{c_\mu}$.

*Part 4*: To exploit the observations of Part 3, it is necessary to choose an appropriate orthonormal frame for $\Lambda^+|_p$. To this end, fix $r \in (0, r_0]$ for the moment, and define the $3 \times 3$ symmetric, non-negative definite matrix $\mathbb{T}(r)$ by declaring its components $\{\mathbb{T}_{ab}(r)\}_{a,b \in \{1,2,3\}}$ to be:

$$\mathbb{T}_{ab}(r) = \frac{1}{r^3} \int_{\partial B_r} \langle a_a a_b \rangle .$$

(4.9)

The trace of $\mathbb{T}$ is $e^{-2\mathfrak{o}} K^2(r)$ which is $K^2(r)$ up to a factor that differs from 1 by $c_0 r^2$. Therefore, $K^2(r)$ is no smaller than $\mathbb{T}$'s largest eigenvalue (up to this same $1 + \mathcal{O}(r^2)$ factor). Also, $\mathbb{T}$'s largest eigenvalue is no smaller than $\frac{1}{3} K^2(r)$ (up to the same factor) because $\mathbb{T}$ is non-negative definite.

The unit length eigenvectors of $\mathbb{T}(r)$ are orthogonal and span $\Lambda^+|_p$ because $\mathbb{T}$ is symmetric. Therefore, there is an oriented, orthonormal basis for $\Lambda^+|_p$ composed of these eigenvectors. This is the basis that will be used to exploit (4.8). Note in this regard that if $\{\omega^c\}_{c=1,2,3}$ are the unit length eigenvectors of $\mathbb{T}(r)$, then $\mathbb{T}_{ab}(r) = 0$ if $a \neq b$. Moreover, the diagonal elements of $\mathbb{T}(r)$ are the numbers $\{K_c^2(r)\}_{c=1,2,3}$; so these are the eigenvalues of $\mathbb{T}(r)$. The convention used subsequently labels the frame $\{\omega^c\}_{c=1,2,3}$ so that the eigenvalues of $\mathbb{T}(r)$ obey $K_1^2 \leq K_2^2 \leq K_3^2$.

An important point for what follows is that $\mathbb{T}$ varies smoothly with $r$. In fact, Stoke's theorem and (3.16) can be used to bound the $r$-derivative of $\mathbb{T}$ according to

$$|\tfrac{d}{dr} \mathbb{T}| \leq c_0 (\tfrac{N}{r} + r) |\mathbb{T}|.$$

(4.10)

This derivative bound has the following implications: Let $\{\omega^c\}_{c=1,2,3}$ denote the basis of eigenvectors of $\mathbb{T}(r)$. If $N(r) \le 1$ and if $s \in [\frac{1}{2}r, r]$, then

$$|\mathbb{T}_{ab}(s) - \delta_{ab} K_a^2(r)| \le c_0 (N(r) + r^2) K^2(r) .$$

(4.11)

Indeed, this follows from (4.10) because of the bound $|\mathbb{T}| \le c_0 K^2$, and because $K$ is increasing, and because $N$ obeys (4.4) (invoke the version with what is denoted by s in (4.4) replaced by $\frac{1}{2}r$.)

*Part 5*: Now suppose that $\varepsilon \in (0, 1]$ and that $N(r) \le \varepsilon$. Suppose in addition that $c \in \{1, 2, 3\}$, that $z > 10$ and that $K_c(r) \ge \frac{1}{zc} K(r)$. (In the applications to come, the number $z$ will be chosen so that $z \le c_\mu$.) This bound is certainly obeyed by $K_3(r)$ because $K_3(r)$ is the largest of the numbers $\{K_c(r)\}_{c=1,2,3}$, and because $K^2(r) = e^{-2\partial} \sum_{c=1,2,3} K_c^2(r)$. Anyway, granted that $K_c(r) \ge \frac{1}{zc} K(r)$, then by virtue of (4.11),

$$|K_c(\cdot) - K_c(r)| \le c_0 zc (\varepsilon + r^2) K(r) \quad on \ [\tfrac{1}{2}r, r];$$

(4.12)

and therefore

$$K_c(\cdot) \ge \frac{1}{2zc} K(r) \quad on \ [\tfrac{1}{2}r, r]$$

(4.13)

if $\varepsilon < c_0^{-1} \frac{1}{z^2 c^2}$ and if $r \le c_0^{-1} \frac{1}{zc}$. Hence, if the matrix $\mathbb{T}(r)$ has two eigenvalue that are greater than $\frac{1}{zc} K(r)$, then the matrix $\mathbb{T}(s)$ for $s \in [\frac{1}{2}r, r]$ has at least two eigenvalues that are greater than $\frac{1}{2zc} K(r)$. Keep in mind in this regard that $K(r) \ge K(s)$ when $s \le r$. Also,

$$K(s) \ge (1 - c_0 \varepsilon) K(r) \quad for \ s \in [\tfrac{1}{2}r, r]$$

(4.14)

when $N(r) \le \varepsilon$ by virtue of (4.4) and (3.7).

The take-away here is that $\mathbb{T}(s)$ has two or more relatively large eigenvalues if $\mathbb{T}(r)$ does.

*Part 6*: Fix $s \in [\frac{1}{2}r, r]$ and let $\lambda_1 \le \lambda_2 \le \lambda_3$ denote the eigenvalues for $\mathbb{T}(s)$. Fix for the moment a number $\delta \in (0, \frac{1}{1000})$. Since

$$2\lambda_3(\lambda_1 + \lambda_2) + 2\lambda_1\lambda_2 = (\text{trace}(\mathbb{T}(s)))^2 - \text{trace}(\mathbb{T}(s)^2) ,$$

(4.15)

a bound saying that $(\operatorname{trace}(\mathbb{T}))^2 - \operatorname{trace}(\mathbb{T}^2) \le \delta \kappa^4(s)$ implies that $\mathbb{T}(s)$ has but one eigenvalue greater than $10\,\delta \kappa^2(s)$. (Remember that $\lambda_3 \ge \frac{1}{3}(1-c_0 r^2)\kappa^2(s)$ and that the eigenvalues of $\mathbb{T}(s)$ are non-negative.)

To see about a bound for the right hand side of (4.15), consider first the square of the trace of $\mathbb{T}$: Keeping in mind that $\operatorname{trace}(\mathbb{T}(s)) = \kappa^2(s)$, and keeping in mind the definition of $\kappa^2(s)$ from Section 3a, one has (using the Cauchy-Schwarz inequality):

$$\frac{1}{2\pi^2}(\operatorname{trace}(\mathbb{T}(s)))^2 \le (1+c_0 r^2)\frac{1}{s^3}\int_{\partial B_s}|a|^4 \ .$$

(4.16)

Noting from the top bullet of Proposition 3.1 that $|a|^4$ on $\partial B_s$ obeys $|a|^4 \le c_\mu \kappa^4(r)$, and noting (4.14), this bound (when $\varepsilon \le c_0^{-1}$) implies in turn the bound

$$\frac{1}{2\pi^2}(\operatorname{trace}(\mathbb{T}(s)))^2 \le \frac{1}{s^3}\int_{\partial B_s}|a|^4 \ + \ c_0 r^2 \kappa^4(s) \ .$$

(4.17)

To see about the trace of $\mathbb{T}^2$, introduce a $3 \times 3$ matrix valued function on $\partial B_s$ to be denoted by $\mathcal{K}$ by using the following rule for its components:

$$\mathcal{T}'_{ab} = \langle a_a a_b \rangle - \frac{1}{2\pi^2}\mathbb{T}_{ab} \ .$$

(4.18)

By virtue of this definition,

$$\frac{1}{2\pi^2}\operatorname{trace}(\mathbb{T}(s)^2) = \frac{1}{s^3}\int_{\partial B_s}|\mathcal{T}'|^2 \ - \ \frac{1}{s^3}\int_{\partial B_s}|\mathcal{T}' - \frac{1}{2\pi^2}\mathbb{T}(s)|^2 + \mathfrak{e} \ ,$$

(4.19)

with $\mathfrak{e}$ having norm bounded by $c_0 r^2 \kappa^4(r)$. This term $\mathfrak{e}$ accounts for the fact that the area of $\partial B_s$ is not necessarily exactly $2\pi^2 s^3$. If it were exactly $2\pi^2 s^3$, then $\mathfrak{e}$ would be zero. Since $|\mathcal{T}'|^2 = \langle a_a a_b \rangle \langle a_a a_b \rangle$, the identity in (4.19) and the inequality in (4.17) give the bound

$$\frac{1}{2\pi^2}((\operatorname{trace}(\mathbb{T}))^2 - \operatorname{trace}(\mathbb{T}^2))|_s \le \frac{1}{s^3}\int_{\partial B_s}(|a|^4 - \langle a_a a_b \rangle \langle a_a a_b \rangle) + \frac{1}{s^3}\int_{\partial B_s}|\mathcal{T}' - \frac{1}{2\pi^2}\mathbb{T}(s)|^2 + c_0 r^2 \kappa^4(r).$$

(4.20)

Granted this last inequality, the task now is to obtain suitable bounds for the two integrals on its right hand side from the assumption that $N$ is small.

*Part 7*: The bound $N(r) < \varepsilon$ has the following additional implication: There is a subset of $[(1-2\mu)r, (1-\mu)r]$ of measure greater than $\frac{99}{100}\mu r$ that is characterized as follows: A point s is in this subset when

$$\frac{1}{s} \int_{\partial B_s} (|\nabla_A a|^2 + r^2 \|[a;a]\|^2) \leq 10^6 \frac{1}{\mu} \varepsilon K^2(r) \ .$$

(4.21)

This inequality will now be used to bound the right hand side of (4.20).

In the first place, (4.21) implies directly

$$\frac{1}{r^3} \int_{\partial B_s} \|[a;a]\|^2 \leq \frac{10^4}{m^2 \mu} \varepsilon K^4(r)$$

(4.22)

with $m$ denoting here $r\, r\, K(r)$. Since $\frac{1}{4}\|[a;a]\|^2 = |a|^4 - \langle a_a a_b \rangle \langle a_a a_b \rangle$, and since $s \geq \frac{1}{2} r$, this inequality implies in turn that

$$\frac{1}{s^3} \int_{\partial B_s} (|a|^4 - \langle a_a a_b \rangle \langle a_a a_b \rangle) \leq c_0 \frac{1}{m^2 \mu} \varepsilon K^4(r) \ .$$

(4.23)

As explained momentarily, this last inequality is suitable for dealing with the left most integral on the right hand side of (4.20).

*Part 8*: To see about the $|\mathcal{T}' - \frac{1}{2\pi^2}\mathbb{T}(s)|^2$ integral in (4.20), let $t_{ab}$ denote the average of the function $\mathcal{T}'_{ab} - \frac{1}{2\pi^2}\mathbb{T}_{ab}(s)$ over $\partial B_s$. This number is zero when the area of $\partial B_s$ is exactly $2\pi^2 s^3$. In general, $|t_{ab}| \leq c_0 r^2$ because the area of $\partial B_s$ differs from $2\pi^2 s^3$ by at most $c_0 s^5$. Meanwhile,

$$\int_{\partial B_s} |\mathcal{T}' - \frac{1}{2\pi^2}\mathbb{T} - t|^2 \leq c_0 s^2 \int_{\partial B_s} |d\mathcal{T}'|^2$$

(4.24)

because the constant functions span the kernel of the operator $d^\dagger d$ on $\partial B_s$. Now, $|d\mathcal{T}'|$ is observedly less than $c_0 |a| |\nabla_A a|$; and $|a|$ is less than $c_\mu K(r)$ (by virtue of the top bullet in Proposition 3.1 since $s \leq (1-\mu)r$). Therefore, (4.24) implies in turn that

$$\int_{\partial B_s} |\mathcal{T}' - \frac{1}{2\pi^2}\mathbb{T} - t|^2 \leq c_0 s^2 K^2(r) \int_{\partial B_s} |\nabla_A a|^2 \ .$$

(4.25)

Noting (4.21), and noting again that $|t| \leq c_0 r^2 K^2(r)$, this last bound leads to:

$$\frac{1}{s^3} \int_{\partial B_s} |\mathcal{T}' - \frac{1}{2\pi^2}\mathbb{T}(s)|^2 \leq c_0(\frac{1}{\mu}\varepsilon + r^2)\kappa^4(r) ,$$

(4.26)

which is the desired bound for the $|\mathcal{T}' - \frac{1}{2\pi^2}\mathbb{T}(s)|^2$ integral in (4.20).

*Part* 9: Putting together the conclusions of Parts 5-8 (and in particular, (4.15), (4.20), (4.23) and (4.26)) gives a bound for $\lambda_3(\lambda_1 + \lambda_2)$ that reads

$$\lambda_3(\lambda_1 + \lambda_2) \leq c_\mu (\varepsilon + r^2)\kappa^4(r) .$$

(4.27)

Since $\lambda_3 \geq \frac{1}{3}(1 - c_0 r^2)\kappa^2(s)$ which is greater than $\frac{1}{6}\kappa^2(r)$ when $\varepsilon \leq c_0^{-1}$ and $r \leq c_0^{-1}$, this bound implies that

$$\lambda_1 + \lambda_2 \leq c_\mu(\varepsilon + r^2)\kappa^2(r) .$$

(4.28)

Hence, both $\lambda_1$ and $\lambda_2$ are less than $\frac{1}{100 z^2 c^2}\kappa^2(r)$ if $\varepsilon < c_\mu^{-1}\frac{1}{z^2 c^2}$ and if $r \leq c_\mu^{-1}\frac{1}{zc}$. This implies in turn that $\mathbb{T}(r)$ has only 1 eigenvalue greater than $\frac{1}{zc}\kappa^2(r)$. Therefore, the right hand side of (4.8) will be much less than $\frac{1}{c^2}$ if $z \geq c_\mu$ and if $\varepsilon \leq c_\mu^{-1}\frac{1}{c^2}$ and $r \leq c_\mu^{-1}\frac{1}{c}$. And, this implies that $r_{c_\mu} > r$.

## 5. When $r \leq r_{c_\mu}$

Supposing that $\mu \in (0, \frac{1}{100}]$ has been chosen and also $c \geq c_\mu$, let p denote a point in X where $r_{c_\mu} > 0$. (Keep in mind that this is the case if $a \neq 0$ at p). The propositions in this section give a priori bounds on the behavior of the components of $a$ and those of $[a; a]$ and $[F_A^-, a]$ on balls of radius less than the value of $r_{c_\mu}$ at the point p.

### a) The variation of $|a|$

To set the stage, fix an orthonormal frame for $\Lambda^+$ at p and, as before, parallel transport it along the geodesic arcs from p using the Levi-Civita connection to define an orthonormal frame for $\Lambda^+$ on the whole of the radius $r_0$ ball centered at p. The proposition that follows makes an assertion to the effect that $|a|$ differs little from $\kappa(r)$ on $B_{(1-\mu)t}$ when $r \leq r_{c_\mu}$. It has the following implication: If $r_{c_\mu} > 0$ for a given point p, then $a \neq 0$ at p. (As noted at the start of Section 4, the converse assertion, that $r_{c_\mu} > 0$ if $a \neq 0$ at p, follows from the definition of $r_{c_\mu}$ because A and $a$ are smooth.)

**Proposition 5.1**: *There exists $\kappa > 1$, and given $\mu \in (0, \frac{1}{100}]$, there exists $\kappa_\mu > \kappa$ that depends only on $\mu$; and these numbers have the following significance: Fix $r > 1$ and suppose that $(A, a)$ is a solution to (1.2). Fix a point $p \in X$ (in $\underline{X}$ if $X$ is not compact) and fix $c \geq \kappa_\mu$ to define $r_{c_\mu}$. Assume that $p$ is such that $r_{c_\mu} > 0$. If $r$ is positive, but at most the smaller of $r_{c_\mu}$ and $\frac{1}{c}$, then*

- $||a| - \frac{1}{\sqrt{2\pi}} K(r)| \leq \kappa_\mu \frac{1}{c} K(r)$ *on* $B_{(1-\mu)r_c}$.

- *Write $a = a_c \omega^c$ using the chosen orthonormal frame for $\Lambda^+$ near $p$. Each component $a_e \in \{a_c\}_{c=1,2,3}$ obeys* $\max_{B_{(1-\mu)r}} |a_e|^2 - \min_{B_{(1-\mu)r}} |a_e|^2 \leq \kappa_\mu \frac{1}{c} K^2(r)$.

***Proof of Proposition 5.1***: To start, note that $N(r) \leq \frac{1}{c^2}$ because $N(r)$ is less than the $q = p$ version of the integral on the left hand side of (4.1). Consider now the top bullet of the proposition: The inequality $|a| \leq \frac{1}{\sqrt{2\pi}} (1 + c_\mu \frac{1}{c}) K(r)$ follows from the top bullet of Proposition 3.1 (given what was just said about $N(r)$.) The other half of the inequality, $|a| \geq \frac{1}{\sqrt{2\pi}} (1 - c_\mu \frac{1}{c}) K(r)$, follows from (3.13) given that $N(r) \leq \frac{1}{c^2}$, given that $r \leq r_{c_\mu}$ and given the top bullet in (1.10).

To start the proof of the second bullet, return to (3.23). With $r$ fixed in advance, the number $K_c(r)$ is now viewed as being just a number, which is to say that $dK_c(r) = 0$. With this understood, then (3.23) and the second bullet of (3.11) lead to an analog of (3.12) that says

$$\frac{1}{2} (|a_c|^2(q) - \frac{1}{2\pi^2} K_e^2(r)) + \int_{B_r} \chi_\mu G_q (|\nabla_A a_c|^2 + r^2 |[a; a_c]|^2)$$
$$= \int_{B_{(1-\mu/4)r}} L_{\mu, q} (|a_c|^2 - K_c^2(r)) + \int_{B_{(1-\mu)r}} \chi_\mu G_q M .$$

(5.1)

There are three integrals in this identity, two on the right and one on the left. The claim by the proposition follows by showing that all of three have small absolute value if $c$ is greater than $c_\mu$; and if $r$ is then less than $r_{c_\mu}$ and $\frac{1}{c}$, and $N(r) < \frac{1}{c^2}$.

Consider first the integral of $\chi_\mu G_q M$: As explained in the proof of Proposition 3.2, the absolute value of this integral is at most $c_0 r^2 K^2(r)$. Consider next the integral with $L_{\mu, q}(|a_c|^2 - K_c^2(r))$. Arguments exactly like those used in Part 2 of the proof of Proposition 3.1 in [T3] bound the absolute value of this integral by $c_\mu (r^2 + \sqrt{N_c(r)}) K_c^2$. Since $N_c K_c^2 \leq N K^2$ (look at the definitions of $N_c$ and $N$), this in turn implies a bound by $c_\mu (r^2 + \sqrt{N}) K K_c$ for the norm of the $L_{\mu, q}$ integral on the right hand side of (5.1).

The final integral to consider is the integral on the left hand side of (5.1). This one is positive, and it is no larger than the result of replacing $a_c$ by $a$ in both occurrences. Therefore, because $r \leq r_{c_\mu}$, this is no greater than $\frac{1}{c^2} K^2(r)$.

Putting these bounds together leads from (5.1) to the bound

$$||a_c|^2 - K_c(r)^2| \leq c_\mu (\tfrac{1}{c} + r^2 + \sqrt{N(r)}) K^2(r)$$

(5.2)

at each point in $B_{(1-\mu)r}$; and this leads immediately to what is asserted by the second bullet of Proposition 5.1.

### b) Writing $a$ as $\nu\sigma_\ddagger + \mathfrak{a}$

Fix a number $\mu \in (0, \tfrac{1}{100}]$ and let $\kappa_*$ denote the largest of the versions of the number $\kappa_\mu$ that appear in Propositions 3.1, 3.2 and 5.1. Fix $p \in X$ (in $\underline{X}$ if X is not compact) and then fix $c > 100\kappa_*$ so as to define the number $r_{c_\mu}$. The upcoming Proposition 5.2 talks about a certain decomposition of $a$ on balls of radius $r < r_{c_\mu}$ subject to the constraint that $rrK(r) \geq 1$.

**Proposition 5.2**: *There exists $\kappa > 1$, and given $\mu \in (0, \tfrac{1}{100}]$, there exists $\kappa_\mu > \kappa$ that depends only on $\mu$; and these numbers have the following significance: Fix $r > 1$ and suppose that $(A, a)$ is a solution to (1.2). Fix a point $p \in X$ (in $\underline{X}$ if X is not compact) and fix $c \geq \kappa_\mu$ to define $r_{c_\mu}$. Assume that p is such that $r_{c_\mu} > 0$. Fix r so that $rrK(r) > \kappa_\mu$, but assume that it is less than both $r_{c_\mu}$ and $\tfrac{1}{c}$. Under these circumstances, a on $B_{(1-\mu)r}$ can be written as $a = \nu\sigma_\ddagger + \mathfrak{a}$ where $\nu$ is a section of $\Lambda^+$ and $\sigma_\ddagger$ is a unit length section of $\mathfrak{G}$ and $\mathfrak{a}$ is a section of $\Lambda^+ \otimes \mathfrak{G}$; and where $\nu$, $\sigma_\ddagger$, and $\mathfrak{a}$ have the properties listed below*
- *The section $\mathfrak{a}$ obeys $\langle \sigma_\ddagger \mathfrak{a} \rangle = 0$ and $\langle \nu, \mathfrak{a} \rangle = 0$*
- *The covariant derivative of $\sigma_\ddagger$ obeys the pointwise bound $|\nabla_A \sigma_\ddagger| \leq \kappa K^{-1}(r) |\nabla_A \mathfrak{a}|$.*
- $d\nu = \langle \nabla_A \sigma_\ddagger \wedge \mathfrak{a} \rangle$
- *The norm of $|\mathfrak{a}|$ on $B_{(1-2\mu)r}$ is no greater than $|\mathfrak{a}| \leq \kappa_\mu \tfrac{1}{c} e^{-\sqrt{rrK(r)}/\kappa_\mu} K(r)$.*
- *If $q \in B_{(1-2\mu)r}$, the integral of $\tfrac{1}{\text{dist}(q,\cdot)^2} |\nabla_A \mathfrak{a}|^2$ on $B_{(1-\mu)r}$ is at most $\kappa_\mu \tfrac{1}{c^2} e^{-\sqrt{rrK(r)}/\kappa_\mu} K^2(r)$.*
- *If $q \in B_{(1-2\mu)r}$, the integral of $\tfrac{1}{\text{dist}(q,\cdot)^2} |\nabla \nu|^2$ on $B_{(1-\mu)r}$ is at most $\kappa_\mu \tfrac{1}{c^2} K^2(r)$.*

This proposition is proved momentarily. Note in the meantime that the first bullet of Proposition 5.2 and Proposition 3.1 lead to the a bound for the norm of $[a; a]$ on $B_{(1-2\mu)r}$:

$$\|[a;a]\| \leq c_\mu \tfrac{1}{c} \, e^{-\sqrt{r\tau K(r)}/\kappa} \, K^2(r) \,.$$

(5.3)

Indeed, this bound follows because $[a;a] = v\,[\sigma, \mathfrak{a}] + [\mathfrak{a};\mathfrak{a}]$ and because $|v| \leq |a|$ (which follows from the orthogonality assertions in the first bullet of Proposition 5.2.)

*Proof of Proposition 5.2*: The proof of Proposition 5.2 has seven parts.

    *Part 1*: To set notation: Any given unit length element in $\Lambda^+|_p$ can be parallel transported along the radial geodesics from p by the Levi-Civita connection to define a corresponding section of $\Lambda^+$ over $B_{r_0}$. If the chosen element in $\Lambda^+|_p$ is denoted by $\omega$, then the corresponding section is also denoted below by $\omega$. Keep in mind that the metric pairing between $\omega$ and $a$ (which is denoted by $\langle\omega,a\rangle$) defines a section of $\mathfrak{G}$ over $B_{r_0}$.

    A convention used in what follows assumes that the constant $c$ is large enough so as to invoke Proposition 5.1. More to the point, the number $c$ should be large enough (greater than $c_\mu$ is sufficient) to guarantee the following:

- $\left| |a| - \tfrac{K(r)}{\sqrt{2\pi}} \right| < \tfrac{1}{1000}$ on $B_{(1-\mu)r}$.
- If $\omega \in \Lambda^+|_p$, then $\max_{B_{(1-\mu)r}} |\langle\omega,a\rangle| - \min_{B_{(1-\mu)r}} |\langle\omega,a\rangle| < \tfrac{1}{1000} K(r)$ on $B_{(1-\mu)r}$.

(5.4)

A positive definite, Hermitian endomorphism of $\mathfrak{G}$ over $B_{(1-\mu)r}$ is defined in the upcoming (5.5) using $a_{\dot{+}}$ and the Riemannian metric. This endomorphism is denoted in what follows by $Y$. The definition of $Y$ uses a chosen orthonormal frame for $\Lambda^+$ on $B_{(1-\mu)r}$ and writes the components of $a$ with respect to this frame as $\{a_c\}_{c\in\{1,2,3\}}$. These are sections over $B_{(1-\mu)r}$ of the bundle $\mathfrak{G}$. With this notation understood, define $Y$ by the rule

$$\sigma \to Y(\sigma) = \tfrac{1}{K^2(r)} a_c \langle a_c \sigma \rangle \,,$$

(5.5)

with it understood (as always) that there is an implicit summation over the repeated indices. Note that $Y$ does not depend on the chosen orthonormal frame.

    The following two assertions are respective consequences of the first and second bullets of (5.4) when $c > c_\mu$:

- *The endomorphism $Y$ has at least one eigenvalue that is greater than $\tfrac{1}{80}$ at each point of $B_{(1-\mu)r}$.*
- *Suppose that $Y$ has only one eigenvalue greater than the minimum of $c^{-2}$ and $\tfrac{1}{750}$ at some point in $B_{(1-\mu)r}$ and that the corresponding eigenspace has multiplicity one.*

*Then there is only one eigenvector with eigenvalue greater than the minimum of $c_\mu c^{-1}$ and $\frac{1}{500}$ at all points of $B_{(1-\mu)r}$ and the eigenspace has multiplicity one also.*

(5.6)

Conditions that guarantee the existence of only one $\mathcal{O}(1)$ eigenvalue of $Y$ at all points in $B_{(1-\mu)r}$ are given in following lemma.

**Lemma 5.3**: *There exists $\kappa > 1$, and given $\mu \in (0, \frac{1}{100}]$, there exists $\kappa_\mu > \kappa$ with the following significance: Fix $r > 1$ and a pair $(A, a)$ that obeys (1.2). Fix $c > \kappa_\mu$ and given $p \in X$ (in $\underline{X}$ if X is not compact) where $a \neq 0$, define $r_{c_\mu}$. Suppose that $r$ is less than $r_{c_\mu}$ and $\frac{1}{c}$. If $rrK(r) \geq \kappa$, then $Y$ has only one eigenvalue greater than $\frac{1}{c}$ at each point of $B_{(1-\mu)r}$ and its eigenspace has dimension one.*

**Proof of Lemma 5.3**: Since $r \leq r_c$,

$$\frac{1}{r^2} r^2 \int_{B_r} \|[a;a]\|^2 \leq \frac{1}{c^2} K^2(r)$$

(5.7)

because the alternative would violate the bound in (4.1) for the case $q = p$. This bound can be rewritten to say that:

$$\frac{1}{r^4} \int_{B_r} \|[a;a]\|^2 \leq \frac{1}{c^2 m^2} K^4(r) ,$$

(5.8)

with $m$ denoting here $rrK(r)$. The bound in (5.8) requires a point (actually, lots of points) in $B_{r/2}$ where $\|[a;a]\| \leq c_0 \frac{1}{cm} K^2(r)$; and this can occur only in the event that $Y$ at this point has only one eigenvalue greater than $c_0 \frac{1}{cm} K^2(r)$. Therefore, if $m$ (which is $rrK$) is greater than $c_0$, then there is a point in $B_{(1-\mu)r}$ that can be used to invoke the second bullet of (5.6); and this second bullet leads directly to the conclusions of the lemma.

*Part 2*: Suppose henceforth that $Y$ has just one eigenvalue that is greater than $\frac{1}{500}$ at all points of $B_{(1-\mu)r}$ and that the corresponding eigenspace is 1-dimensional. There is in this case a section of $\mathfrak{G}$ over $B_{(1-\mu)r}$ whose restriction to any given point is an eigenvector of $Y$ with the largest eigenvalue. This eigensection is denoted in what follows by $\sigma_\ddagger$. Use $\lambda_\ddagger$ to denote the corresponding eigenvalue. It follows from Proposition 5.1 that $\lambda_\ddagger$ differs from $\frac{1}{2\pi}$ by at most $c_\mu c^{-1}$, this difference being less than $c^{-1/2}$ when $c \geq c_\mu$.

The section $a$ is written over $B_{(1-\mu)r}$ using $\sigma_\ddagger$ as

$$a = v\sigma_{\ddagger} + \mathfrak{a}$$

(5.9)

with $v$ being a section of and $\Lambda^+$ on $B_{(1-\mu)r}$, and with $\mathfrak{a}$ being a section of the bundle $\Lambda^+ \otimes \mathfrak{G}$ over $B_{(1-\mu)r}$ that obeys $\langle \sigma_{\ddagger} \mathfrak{a} \rangle = 0$ at each point. Moreover,

- $\lambda_{\ddagger} = |v|^2$.
- $|\mathfrak{a}| \leq c_{\mu} \frac{1}{c^{1/2}} K(r)$.
- $\langle v, \mathfrak{a} \rangle = 0$
- $\|[a;a]\|^2 = 4|v|^2|\mathfrak{a}|^2 + |\mathfrak{a} \wedge \mathfrak{a}|^2$.
- $\left| |v| - \frac{K(r)}{\sqrt{2}\pi} \right| \leq c_{\mu} \frac{1}{c^{1/2}} K(r)$
- Let $v_p$ denote the section of $\Lambda^+$ on $B_{(1-\mu)r}$ that is defined by the parallel transport of $v|_p$ along the radial geodesics from p. Then $|v - v_p| \leq c_{\mu} \frac{1}{c^{1/2}} K(r)$ on $B_{(1-\mu)r}$.

(5.10)

The first three bullets of (5.10) are proved by writing $Y$ using (5.9) as the endomorphism

$$\sigma \to |v|^2 \sigma_{\ddagger} \langle \sigma_{\ddagger} \sigma \rangle + \sigma_{\ddagger} v_{\alpha} \langle \mathfrak{a}_{\alpha} \sigma \rangle + \mathfrak{a}_{\alpha} v_{\alpha} \langle \sigma_{\ddagger} \sigma \rangle + \mathfrak{a}_{\alpha} \langle \mathfrak{a}_{\alpha} \sigma \rangle.$$

(5.11)

The first and third bullets of (5.10) follow from the $\sigma = \sigma_{\ddagger}$ version of (5.11) and the second bullet follows from the second bullet of (5.6) by taking $\sigma$ in (5.11) to be pointwise orthogonal to $\sigma_{\ddagger}$. The fourth bullet in (5.11) can be seen by first writing $[a;a]$ using (5.10) and then invoking the third bullet of (5.11) and the fact that $\langle \mathfrak{a}\sigma_{\ddagger} \rangle = 0$. Note in this regard that these two fact imply that $\|[v\sigma_{\ddagger}; \mathfrak{a}]\| = 2|v||\mathfrak{a}|$; and the fact that $\langle \mathfrak{a}\sigma_{\ddagger} \rangle = 0$ implies by itself that $[\mathfrak{a};\mathfrak{a}]$ is everywhere the tensor product of a section of $\Lambda^+$ and $\sigma_{\ddagger}$. The fifth and sixth bullets follow from Proposition 5.1.

*Part 3*: Supposing that $c \geq c_0$, then the derivatives of $v$, $\sigma_{\ddagger}$ and $\mathfrak{a}$ obey

- $|\nabla v| + |\nabla_A \mathfrak{a}| \leq c_0 |\nabla_A a|$.
- $\nabla_A \sigma_{\ddagger} = -\frac{1}{4} \frac{1}{|v|^2} [\sigma_{\ddagger}, [\sigma_{\ddagger}, *(v \wedge *d_A \mathfrak{a})]]$
- $dv = \langle d_A \sigma_{\ddagger} \wedge \mathfrak{a} \rangle$.

(5.12)

The first bullet is proved by differentiating both sides of (5.9) and then using the fact that $\langle \sigma_{\ddagger} \nabla_A \sigma_{\ddagger} \rangle = 0$ and $\langle \sigma_{\ddagger} \mathfrak{a} \rangle$ both vanish. The second second and third bullets of (5.12) follow by writing the equation $d_A a = 0$ using (5.9) as follows:

$$d\nu\sigma_{\ddagger} + \nu \wedge d_A\sigma_{\ddagger} + d_A\mathfrak{a} = 0 \, .$$

(5.13)

The second bullet of (5.12) follows from this identity using three observations: The first is that $*(\nu \wedge *(\nu \wedge \alpha)) = -|\nu|^2 \alpha$ for any 1-form $\alpha$ (this is because $\nu$ is self-dual.) The second is that $[\sigma_{\ddagger}, [\sigma_{\ddagger}, \alpha]] = -4\alpha$ in the case where $\alpha$ is a section of $\mathfrak{G}$ obeying $\langle \sigma_{\ddagger} \alpha \rangle = 0$. The third is that $\langle \sigma_{\ddagger} d_A \sigma_{\ddagger} \rangle = 0$ (because $|\sigma_{\ddagger}| = 1$). The third bullet in (5.12) is obtained from (5.13) by contracting with $\sigma_{\ddagger}$ (keep in mind that $\langle \sigma_{\ddagger} d_A \mathfrak{a} \rangle = -\langle d_A \sigma_{\ddagger} \wedge \mathfrak{a} \rangle$).

The fact that $\langle \sigma_{\ddagger} \mathfrak{a} \rangle = 0$, and the first three bullets of (5.10) with the second and third bullets of (5.12) lead directly to the assertions of Proposition 5.2's first three bullets.

*Part 4*: This part of the proof (and the next two parts) prove the assertion that is made by the fourth bullet of Proposition 5.2. This task begins with (2.1). In particular, (2.1) leads (after significant algebraic manipulations) to a second order equation for $\mathfrak{a}$ (when $c \geq c_\mu$ which makes $|\mathfrak{a}| \leq c_0^{-1}$) that has the schematic form

$$\nabla_A^{\dagger} \nabla_A \mathfrak{a} + 4 r^2 |\nu|^2 \mathfrak{a} + \mathcal{R} = 0 \, ,$$

(5.14)

with $|\mathcal{R}| \leq c_\mu (1 + r^2 |\mathfrak{a}|^2) |\mathfrak{a}|^2 + c_\mu \frac{1}{K(r)} |\nabla_A a| |\nabla_A \mathfrak{a}|$. The derivation of (5.14) is left to the reader with the hint that the identities in $\langle \sigma_{\ddagger} \mathfrak{a} \rangle = 0$ and $\langle \nu, \mathfrak{a} \rangle = 0$ and the others in (5.10) and those in (5.12) all play a significant role. In particular, they lead to formulas that say:

- $\nabla^{\dagger} \nabla \nu + \mathfrak{R} \nu = \langle \nabla_A^{\dagger} \nabla_A \sigma_{\ddagger}, \mathfrak{a} \rangle + \mathfrak{e}_\nu$ where $|\mathfrak{e}_\nu| \leq c_0(r^2 |\mathfrak{a}|^2 + \frac{1}{K(r)} |\nabla_A \mathfrak{a}|^2 |\nabla_A a|)$ .
- $\nabla_A^{\dagger} \nabla_A \sigma_{\ddagger} = |\nu|^{-2} \langle \nabla^{\dagger} \nabla \nu, \mathfrak{a} \rangle + \mathfrak{e}_\sigma$ where $|\mathfrak{e}_\sigma| \leq c_0(r^2 |\mathfrak{a}|^3 + \frac{1}{K^2(r)} |\nabla_A \mathfrak{a}| |\nabla_A a|)$ .

(5.15)

These identities can be used to write both $\nabla^{\dagger} \nabla \nu$ and $\nabla_A^{\dagger} \nabla_A \sigma_{\ddagger}$ in terms of $\mathfrak{e}_\nu$ and $\mathfrak{e}_\sigma$ when $|\mathfrak{a}| \ll |\nu|$ (which is guaranteed when $c \geq c_\mu$). The point of writing them in this way is that neither $\mathfrak{e}_\nu$ nor $\mathfrak{e}_\sigma$ involve second derivatives.

Taking the inner product of both sides of (5.14) with $\mathfrak{a}$ leads to a differential inequality for $|\mathfrak{a}|^2$ on $B_{(1-\mu)r}$ asserting the following:

$$\tfrac{1}{2} d^{\dagger} d |\mathfrak{a}|^2 + \tfrac{1}{2} |\nabla_A \mathfrak{a}|^2 + r^2 K^2(r) |\mathfrak{a}|^2 \leq c_0 \frac{1}{K^2(r)} |\mathfrak{a}|^2 |\nabla_A a|^2 \, .$$

(5.16)

The derivation of this equation assumes that $c > c_\mu$ (so $|\mathfrak{a}| \leq c_0^{-1}$) and that $rK(r) > c_0$. These constraints are assumed henceforth so to exploit (5.16).

*Part 5*:  With (5.16) in mind, let $\mathcal{G}_q$ denote the Dirichelet Green's function with pole at q for the operator $d^\dagger d \mathcal{G}_q + r^2 K^2(r)$ on $B_r$. This function is smooth and strictly positive on $B_r - q$. In addition, it obeys the following by way of a priori bounds on $B_{(1-\mu)r}$:

- $\mathcal{G}_q \leq c_0 \frac{1}{\mathrm{dist}(q,\cdot)^2} e^{-rK(r)\mathrm{dist}(q,\cdot)/c_\mu}$ .

- $|\nabla \mathcal{G}_q| \leq c_\mu (1 + rrK(r)) \frac{1}{\mathrm{dist}(q,\cdot)^3} e^{-rK(r)\mathrm{dist}(q,\cdot)/c_\mu}$ .

(5.17)

Use the basic cut-off function $\chi$ to construct a non-increasing function (to be denoted by $\chi_\ddagger$) that is equal to 1 where $\mathrm{dist}(p, \cdot) \leq (1 - \frac{3}{2}\mu)r$ and equal to 0 where $\mathrm{dist}(p, \cdot) \geq (1 - \frac{5}{4}\mu)r$. This function can and should be constructed so that the norm of its derivative is bounded by $c_0 \frac{1}{r\mu}$, and that of its Hessian is bounded by $c_0(\frac{1}{r\mu})^2$. Having fixed a point q with $\mathrm{dist}(p,q) \leq (1 - \frac{3}{2}\mu)r$, multiply both sides of (5.16) by $\chi_\ddagger \mathcal{G}_q$ and then integrate the resulting inequality over the ball $B_r$. Integration by parts leads to

$$|\mathfrak{a}|^2(q) + \int_{B_{(1-\mu)r}} \mathcal{G}_q |\nabla_A \mathfrak{a}|^2 \leq c_0 \frac{1}{K^2(r)} \int_{B_{(1-\mu)r}} \mathcal{G}_q |\mathfrak{a}|^2 |\nabla_A \mathfrak{a}|^2 + c_0 \int_{B_r} (-d^\dagger d\chi_\ddagger \mathcal{G}_q + 2\langle d\chi_\ddagger, d\mathcal{G}_q \rangle) |\mathfrak{a}|^2$$

(5.18)

The next part of the proof gets the desired bounds on $|\mathfrak{a}|$ from (5.18).

*Part 6*:  If the distance from q to p is less than $(1 - 2\mu)r$, then the distance from q to the support of $d\chi_\ddagger$ is greater than $\frac{1}{2}\mu r$. This fact with the bounds in (5.17) (and the fact that $|\mathfrak{a}|^2 \leq c_\mu \frac{1}{c} K^2(r)$ from (5.10)) leads directly to a bound for integral in (5.18) with derivatives of $\chi_\ddagger$: Its absolute value is at most $c_\mu \frac{1}{c} e^{-rrK(r)/c_\mu} K^2(r)$.

To bound the integral of $\mathcal{G}_q |\mathfrak{a}|^2 |\nabla_A \mathfrak{a}|^2$, introduce by way of notation L to denote the greatest integer less than the square root of $rrK(r)$; and for $k \in \{0, 1, \ldots, N\}$, let $\rho_k$ denote $(1 - (1 + \frac{k}{L})\mu)r$. Thus, $\rho_0 = (1-\mu)r$ and $\rho_L = (1-2\mu)r$. For each such k, define a number $z_k$ by writing the maximum of $|\mathfrak{a}|$ on the ball of radius $\rho_k$ centered at p as $z_k K(r)$. The number $z_0$ is at most $c_\mu \frac{1}{c^{1/2}}$ by virtue of (5.10).

If q is in the radius $\rho_k$ ball centered at p, then the integral of $\mathcal{G}_q |\mathfrak{a}|^2 |\nabla_A \mathfrak{a}|^2$ that appears in (5.18) should be written as a sum of two parts. The first is the contribution from the ball of radius $\rho_{k-1}$ centered at p, and the second is the complement of this ball. These two parts of the integral are bounded respectively by

$$c_\mu \frac{1}{c^2} z_{k-1}^2 K^4(r) \quad and \quad c_\mu \frac{1}{c^3} e^{-rrK(r)/(Lc_\mu)} K^4(r)$$

(5.19)

These bounds are obtained by invoking the top bullet in (5.17); and by using $z_{k-1} K(r)$ for the maximum of $|a|$ on the radius $\rho_{k-1}$ ball centered at p, and using $c_\mu \frac{1}{c^{1/2}} K(r)$ on the complement of this ball. Note that (4.1) is used for both integrals to bound the integral of $\frac{1}{\text{dist}(q,\cdot)^2} |\nabla_A a|^2$.

The bound in (5.19) and the bound for the right most integral in (5.18) lead to a recursive inequality that bounds $z_k^2$ as follows:

$$z_k^2 \le c_\mu \frac{1}{c^2} z_{k-1}^2 + c_\mu \frac{1}{c^3} e^{-\sqrt{rrK(r)}/c_\mu} .$$

(5.20)

(With regards to the exponential term in (5.20), remember that L is at most the square root of $rrK(r)$.) If $c > c_\mu$, then the factor multiplying $z_{k-1}^2$ in (5.20) is less $\frac{1}{c}$. This understood, then (5.20) can be iterated to see that

$$z_k^2 \le (\tfrac{1}{c})^k + c_\mu \frac{1}{c^3} e^{-\sqrt{rrK(r)}/c_\mu} .$$

(5.21)

The k=L version of (5.21) gives the bound that is asserted for $|a|$ by Proposition 5.2.

*Part* 7: This part of the proof addresses the fifth and six bullets of Proposition 5.2. To this end, note that the six bullet of Proposition 5.2 follows directly from (5.12) and (3.13).

To prove the fifth bullet, use the function $\chi$ to construct yet another bump function, this one denoted by $\chi_\diamond$. This function is non-increasing, equal to 1 on the radius $(1 - \tfrac{31}{16} \mu) r$ ball centered at p and equal to 0 on the complement of the radius $(1 - \tfrac{15}{8} \mu) r$ ball centered at p. It should be constructed so that the norm of its derivative is at most $c_0 \frac{1}{r\mu}$ and so that the norm of its Hessian is at most $c_0 (\frac{1}{r\mu})^2$.

Fix $q \in B_{(1-2\mu)r}$ and reintroduce the Green's function $G_q$ from (1.10). Multiply both sides of (5.14) by $\chi_\diamond G_q$ and then integrate the result over $B_r$. Then integrate by parts to remove $d^\dagger d$ from $|a|^2$. Because the $k = \tfrac{1}{2} L$ version of (5.19) leads to a $c_\mu \tfrac{1}{c} e^{-\sqrt{rrK(r)}/c_\mu} K(r)$ bound for $|a|$ on the radius $(1 - \tfrac{15}{8}\mu) r$ ball centered at p, the resulting integral inequality (after the integration by parts) with (1.10)'s three bullets lead directly to the bound for the integral of $\frac{1}{\text{dist}(q,\cdot)^2} |\nabla_A a|^2$ that is asserted by the fifth bullet of Proposition 5.2.

c) **The second derivatives of** *a*

Fix a number $\mu \in (0, \tfrac{1}{100}]$ and let $\kappa_*$ denote the largest of the versions of the number $\kappa_\mu$ that appear in Propositions 3.1, 3.2 and 5.1 and 5.2. Fix $p \in X$ (in $\underline{X}$ if X is not compact) and then fix $c > 100\kappa_*$ so as to define the number $r_{c\mu}$. Supposing that this is positive, take $r < r_{c\mu}$ and write $a$ on $B_r$ as $\nu\sigma_{\ddagger} + \mathfrak{a}$ in the manner of Proposition 5.2. The upcoming Proposition 5.4 gives a bound for the $L^2$ norm of the second (covariant) derivatives of $\nu$, $\sigma_{\ddagger}$ and $\mathfrak{a}$ on the radius $(1-3\mu)r$ ball centered at p.

To set the notation for the proposition, define a non-negative function F on $[0, r_0]$ by writing the integral of $|F_A|^2$ on $B_r$ as $r^2 F^2(r)$. In the case when X is compact, this function F is bounded by $c_0$ because the square of the $L^2$ norm of $F_A^+$ is bounded by $c_0 r^2$ (by virtue of Lemma 2.1) and the difference between the squares of the $L^2$ norms of $F_A^+$ and $F_A^-$ is bounded by a universal multiple of the first Pontrjagin class of the principal bundle P. (The inequality in (5.3) can be used to bound the square of the $L^2$ norm of $F_A^+$ on a compact X by a multiple of power of $r$ that is less than 2. Supposing this, then the norm of F would be smaller than some negative power of $r$.)

**Proposition 5.4**: *There exists $\kappa > 1$, and given $\mu \in (0, \frac{1}{100}]$, there exists $\kappa_\mu > \kappa$ that depends only on $\mu$; and these numbers having the following significance: Fix $r > 1$ and suppose that $(A, a)$ is a solution to (1.2). Fix a point $p \in X$ (in $\underline{X}$ if $X$ is not compact) and fix $c \geq \kappa_\mu$ to define $r_{c\mu}$. Assume that $p$ is such that $r_{c\mu} > 0$. Fix $r$ so that $r r K(r) > \kappa_\mu$, but assume that it is less than both $r_{c\mu}$ and $\frac{1}{c}$. Write $a$ on $B_r$ as $a = \nu\sigma_{\ddagger} + \mathfrak{a}$ in the manner of Proposition 5.3.*

- $\int_{B_{(1-3\mu)r}} |\nabla_A \nabla_A \mathfrak{a}|^2 \leq \kappa_\mu \frac{1}{c^2} K^2(r)(1 + \frac{F^2(r)}{r^2 K^2(r)}) e^{-\sqrt{r r K(r)}/\kappa_\mu}$.

- $\int_{B_{(1-3\mu)r}} |\nabla_A \nabla_A \sigma_{\ddagger}|^2 \leq \kappa_\mu \frac{1}{c^2} (1 + \frac{F^2(r)}{r^2 K^2(r)}) e^{-\sqrt{r r K(r)}/\kappa_\mu}$

- *Let $\nu_{\ddagger}$ denote a given self-dual, closed 2-form on $B_r$ (it can be 0). Then*

$$\int_{B_{(1-3\mu)r}} |\nabla\nabla(\nu - \nu_{\ddagger})|^2 \leq \kappa_\mu (\int_{B_r} |\nu - \nu_{\ddagger}|^2 + \frac{1}{c^2}(K^2(r) + \frac{F^2(r)}{r^2}) e^{-\sqrt{r r K(r)}/\kappa_\mu}).$$

This proposition has the following by way of an immediate consequence.

**Proposition 5.5**: *The assumptions in Proposition 5.4 imply the following bound for the $L^2$ norm of $[F_A, \sigma_{\ddagger}]$ on $B_{(1-3\mu)r}$:* $\int_{B_{(1-3\mu)r}} |[F_A, \sigma_{\ddagger}]|^2 \leq \kappa_\mu \frac{1}{c^2} (1 + \frac{F^2(r)}{r^2 K^2(r)}) e^{-\sqrt{r r K(r)}/\kappa_\mu}$.

The remainder of this section contains the proofs of these propositions (in reverse order).

***Proof of Proposition 5.5***: The proof assumes that Proposition 5.4 is true. Granted Proposition 5.4, the identity $d_A d_A \sigma_\ddagger = [F_A, \sigma_\ddagger]$ can be used to bound the $L^2$ norm of $[F_A, \sigma_\ddagger]$ by $c_0$ times the $L^2$ norm of $\nabla_A d_A \sigma_\ddagger$; and the second bullet of Proposition 5.4 bounds the latter.

***Proof of Proposition 5.4***: The proof of the proposition has six parts. These parts assume at the outset that $c$ is large enough (greater than $c_\mu$) so that Propositions 5.1 and 5.2 can be invoked. Parts 1-5 prove the first and second assertions, and Part 6 proves the third.

*Part 1*: The second bullet of the proposition follows from the first and third (and (1.9) and Proposition 5.1). To elaborate: By virtue of second bullet in (5.12) (differentiate the identity) the second derivative of $\sigma_\ddagger$ is bounded in turn by

$$|\nabla_A \nabla_A \sigma_\ddagger| \leq c_0 \left( \frac{1}{K^2(r)} |\nabla v| |\nabla_A \mathfrak{a}| + \frac{1}{K(r)} |\nabla_A \nabla_A \mathfrak{a}| \right)$$

(5.22)

(The derivation uses the fifth bullet of (5.10) and it assumes that $c \geq c_\mu$.) Therefore, the square of the $L^2$ norm of $\nabla_A \nabla_A \sigma_\ddagger$ is bounded by $c_0$ times the following:

$$\int_{B_{(1-3\mu)r}} |\nabla_A \nabla_A \sigma_\ddagger|^2 \leq \frac{1}{K^4(r)} \left( \int_{B_{(1-3\mu)r}} |\nabla v|^4 \right)^{1/2} \left( \int_{B_{(1-3\mu)r}} |\nabla_A \mathfrak{a}|^4 \right)^{1/2} + \frac{1}{K^2(r)} \int_{B_{(1-3\mu)r}} |\nabla_A \nabla_A \mathfrak{a}|^2 \right).$$

(5.23)

(The derivation uses the Cauchy-Schwarz inequality.) Given the top bullet in (1.9), and given the $L^2$ bounds in Propositions 5.2, and given those in the first and third bullets of Proposition 5.4, then the bound in (5.23) for the square of the $L^2$ norm of $\nabla_A \nabla_A \sigma_\ddagger$ leads directly to the bound that is asserted by the second bullet of Proposition 5.4. (Note that the inequality in (5.23) is used again to prove the third bullet of Proposition 5.2.)

*Part 2*: Use the standard cut-off function $\chi$ to construct a non-increasing function (to be denoted by $\chi_\#$) that equals 1 on $B_{(1-3\mu)r}$ and equals 0 on the complement of $B_{(1-5\mu/2)r}$. This function can and should be constructed so that the norm of its differential is bounded by $c_0 \frac{1}{r\mu}$. Return to (5.16) and multiply both sides of this inequality by $\chi_\#^4 |\nabla_A \mathfrak{a}|^2$. An integration by parts (to remove one derivative from the factor $d^\dagger d|\mathfrak{a}|^2$) leads to:

$$\int_{B_{(1-2\mu)r}} \chi_\#^4 |\nabla_A \mathfrak{a}|^4 \leq c_\mu \frac{1}{K^2(r)} \int_{B_{(1-2\mu)r}} \chi_\#^4 |\mathfrak{a}|^2 |\nabla_A \mathfrak{a}|^2 |\nabla v|^2 + c_\mu \int_{B_{(1-2\mu)r}} \chi_\#^4 |\mathfrak{a}| |\nabla_A \mathfrak{a}|^2 |\nabla_A \nabla_A \mathfrak{a}|$$
$$+ c_\mu \frac{1}{r} \int_{B_{(1-2\mu)r}} |\mathfrak{a}| |\nabla_A \mathfrak{a}|^3 .$$

(5.24)

By way of an explanation, the left most integral on the right hand side of (5.24) accounts for the term on the right hand side of (5.16) because $|\nabla_A a|^2 \leq c_0(|\nabla v|^2 + |\nabla_A \mathfrak{a}|^2)$ (by virtue of the second bullet in Proposition 5.2) and the $|\nabla_A \mathfrak{a}|^2$ part of this contributes a term that is much less than the left hand side of (5.24) (because of the apriori bound on $|\mathfrak{a}|$). Meanwhile, the middle and the right most integrals in (5.24) account for the integration by parts that takes one derivative from $d^\dagger d|\mathfrak{a}|^2$ and puts it on either a factor $|\nabla_A \mathfrak{a}|$ (the middle integral) or a factor of $\chi_\#$ (the far right hand integral).

*Part 3*: Use the bound from the fourth bullet in Proposition 5.2 on $|\mathfrak{a}|$, and the bound in the fifth bullet of Proposition 5.2 to see that (5.24) leads in turn to:

$$\int_{B_{(1-2\mu)r}} \chi_\#^4 |\nabla_A \mathfrak{a}|^4 \leq c_\mu \frac{1}{c^2} e^{-\sqrt{r r K(r)}/c_\mu} \left( \int_{B_{(1-2\mu)r}} \chi_\#^4 |\nabla v|^4 + K^2(r) \int_{B_{(1-2\mu)r}} \chi_\#^4 |\nabla_A \nabla_A \mathfrak{a}|^2 + K^4(r) \right).$$

(5.25)

There are two integrals that appear on the right hand side of (5.25). The left most one is the integral of the fourth power of $\chi_\#|\nabla v|$. To bound this integral, use standard elliptic $L^4$ bounds for the operator d acting on $C^\infty(B_r; \Lambda^+)$ to see that

$$\int_{B_{(1-2\mu)r}} \chi_\#^4 |\nabla v|^4 \leq c_\mu \int_{B_{(1-2\mu)r}} \chi_\#^4 |dv|^4 + c_\mu \frac{1}{r^4} \int_{B_{(1-2\mu)r}} |v|^4 .$$

(5.26)

Keep in mind that $\chi_\#$ has support in $B_{(1-5\mu/2)r}$. (See for example Theorem 6.2.6 in [Mo]; or Theorem 10.3.1 and Corollary 10.3.3 in [Ni].) Given the second and third bullets of Proposition 5.2 (and Proposition 3.1's bound on $|a|$), the inequality in (5.26) leads to

$$\int_{B_{(1-2\mu)r}} \chi_\#^4 |\nabla v|^4 \leq c_\mu \frac{1}{c^4} e^{-\sqrt{r r K(r)}/c_\mu} \int_{B_{(1-2\mu)r}} \chi_\#^4 |\nabla_A \mathfrak{a}|^4 + c_\mu K^4(r) .$$

(5.27)

Using this bound in (5.25) leads to the bound:

$$\int_{B_{(1-2\mu)r}} \chi_\#^4 |\nabla_A \mathfrak{a}|^4 \leq c_\mu \frac{1}{c^2} e^{-\sqrt{r r K(r)}/c_\mu} K^2(r) \left( \int_{B_{(1-2\mu)r}} \chi_\#^4 |\nabla_A \nabla_A \mathfrak{a}|^2 + K^2(r) \right) .$$

(5.28)

Holding on to (5.28), return now to (5.14). Take the square of the norm of both sides, and multiply both sides of the resulting inequality by $\chi_\#^4$. Doing so (and using the third and fourth bullets of Proposition 5.1) leads to the pointwise inequality

$$\chi_\#^4 |\nabla_A^\dagger \nabla_A \mathfrak{a}|^2 \le c_0 (\tfrac{1}{c^4} e^{-\sqrt{r\tau K(r)}/c_\mu} r^4 K^6(r) + \tfrac{1}{K^2(r)} (\chi_\#^4 |\nabla v|^2 |\nabla_A \mathfrak{a}|^2 + \chi_\#^4 |\nabla_A \mathfrak{a}|^4)) .$$

(5.29)

Because of (5.27) and (5.28), integrating both sides of this inequality gives the bound

$$\int_{B_{(1-2\mu)r}} \chi_\#^4 |\nabla_A^\dagger \nabla_A \mathfrak{a}|^2 \le c_\mu \tfrac{1}{c^2} e^{-\sqrt{r\tau K(r)}/c_\mu} K^2(r) + c_\mu \tfrac{1}{c^2} e^{-\sqrt{r\tau K(r)}/c_\mu} \int_{B_{(1-2\mu)r}} \chi_\#^4 |\nabla_A \nabla_A \mathfrak{a}|^2 .$$

(5.30)

The ramifications of this inequality are explored in the next part of the proof.

*Part 4*: Four back-and-forth integration by parts applied to the $|\nabla_A^\dagger \nabla_A \mathfrak{a}|^2$ integral on the left hand side leads to the bound

$$\int_{B_{(1-2\mu)r}} \chi_\#^4 |\nabla_A \nabla_A \mathfrak{a}|^2 \le \int_{B_{(1-2\mu)r}} \chi_\#^4 |\nabla_A^\dagger \nabla_A \mathfrak{a}|^2 - 2 \int_{B_{(1-2\mu)r}} \chi_\#^4 \langle F_{A\alpha\beta} [\nabla_{A\alpha} \mathfrak{a}_c, \nabla_{A\beta} \mathfrak{a}_c] \rangle$$
$$+ c_0 \int_{B_{(1-2\mu)r}} \chi_\#^4 |\mathfrak{a}| |\nabla_A \mathfrak{a}| |d_A^\dagger F_A| + c_0 \tfrac{1}{r} \int_{B_{(1-2\mu)r}} \chi_\#^3 \|[F_A \mathfrak{a}]\| |\nabla_A \mathfrak{a}| + c_0 \tfrac{1}{r^2} \int_{B_{(1-2\mu)r}} |\nabla_A \mathfrak{a}|^2 .$$

(5.31)

This leads (via integration by parts and the Bianchi identity $d_A^\dagger F_A = 2 d_A^\dagger F_A^+$) to

$$\int_{B_{(1-2\mu)r}} \chi_\#^4 |\nabla_A \nabla_A \mathfrak{a}|^2 \le \int_{B_{(1-2\mu)r}} \chi_\#^4 |\nabla_A^\dagger \nabla_A \mathfrak{a}|^2 + c_0 \int_{B_{(1-2\mu)r}} \chi_\#^4 \|[F_{A\alpha\beta}, \mathfrak{a}]\|^2$$
$$+ c_0 \tfrac{1}{r^2} \int_{B_{(1-2\mu)r}} |\nabla_A \mathfrak{a}|^2 + c_0 r^2 \int_{B_{(1-2\mu)r}} |\mathfrak{a}|^2 |\nabla_A \mathfrak{a}|^2 .$$

(5.32)

Plugging in (5.30) to deal with the $|\nabla_A^\dagger \nabla_A \mathfrak{a}|^2$ integral in (5.32) (and assuming that $c > c_\mu$), and using the bounds from Proposition 5.2 to deal with the $|\nabla_A \mathfrak{a}|^2$ and $|\mathfrak{a}|^2 |\nabla_A \mathfrak{a}|^2$ integrals leads next from (5.32) to the bound

$$\int_{B_{(1-2\mu)r}} \chi_\#^4 |\nabla_A \nabla_A \mathfrak{a}|^2 \le c_\mu \tfrac{1}{c^2} e^{-\sqrt{r\tau K(r)}/c_\mu} K^2(r) + c_0 \int_{B_{(1-2\mu)r}} \chi_\#^4 \|[F_{A\alpha\beta}, \mathfrak{a}]\|^2 .$$

(5.33)

To elaborate some, the $|\mathfrak{a}|^2 |\nabla_A \mathfrak{a}|^2$ integral on the right hand side of (5.32) makes a contribution to the right hand side that is no larger than the factor $(\tfrac{1}{r^2} + r^2)$ times the square of the sup norm of $\mathfrak{a}$ times the square of the $L^2$ norm of $|\nabla_A \mathfrak{a}|$. Because of the third and fourth bullets of Proposition 5.2, t, this contribution is no larger than

$$c_\mu \tfrac{1}{c^4} (\tfrac{1}{r^2} + r^2) r^2 e^{-\sqrt{r\tau K(r)}/c_\mu} K^4(r). \text{ This,}$$

(5.34)

What is written in (5.34) is no greater than $\frac{1}{c^2} e^{-\sqrt{r r_K(r)}/c_\mu} K^2(r)$ because of the inequality $x^2 e^{-x} \leq 16 e^{-x/2}$ applied to the case $x = \sqrt{r r_K(r)}/c_\mu$.

*Part 5*: The final stem to proving the first bullet's assertion in Proposition 5.4 bounds the integral $|[F_A, \mathfrak{a}]|$ integral in (5.33) by the product of the square of the sup norm for $|\mathfrak{a}|$ on $B_{(1-2\mu)r}$ (which comes from Proposition 5.2) times the integral of $|F_A|^2$ on $B_{(1-2\mu)r}$. Because the latter is at most $r^2 F^2(r)$, this gives

$$\int_{B_{(1-2\mu)r}} \chi_\#^4 \, |[F_{A\alpha\beta}, \mathfrak{a}]|^2 \leq c_\mu \frac{1}{c^2} e^{-\sqrt{r r_K(r)}/c_\mu} r^2 K^2(r) F^2(r).$$

(5.35)

Using this bound in (5.33) gives the bound that is asserted by Proposition 5.4. (Write $r^2 K^2(r)$ as $\frac{1}{r^2} \times (r^2 r^2 K^2(r))$ and use the $x^2 e^{-x} \leq 16 e^{-x/2}$ bound again with $x = \sqrt{r r_K(r)}/c_\mu$ to replace the factor of $r^2$ with the factor $\frac{1}{r^2}$ at the cost of a larger version of $c_\mu$.)

*Part 6*: To prove the second assertion, note that if $v_\ddagger$ is self-dual and if $dv_\ddagger = 0$, then it obeys $\nabla^\dagger \nabla v_\ddagger + \mathfrak{R} v_\ddagger = 0$. Subtract this equation from the identity depicted by the top bullet in (5.15) to see that

$$|\nabla^\dagger \nabla (v - v_\ddagger)|^2 \leq c_0 |\nabla_A \nabla_A \sigma_\ddagger|^2 |\mathfrak{a}|^2 + c_0 |\nabla_A \sigma_\ddagger|^2 |\nabla_A \mathfrak{a}|^2 + c_0 |v - v_\ddagger|^2 \,.$$

(5.36)

Now multiply both sides of this inequality by $\chi_\#^2$, square both sides and then integrate over $B_{(1-2\mu)r}$. Then, invoke (5.23) and (5.27) and (5.28) and the first assertion in Proposition 5.4 to bound the integral of $\chi_\#^4 |\nabla^\dagger \nabla (v - v_\ddagger)|^2$. A straightforward sequence of back and forth integration by parts leads from the latter bound to the bound that is claimed by the second assertion in Proposition 5.4.

**d) The curvature of A**

Fix $r > 1$ and let $(A, \mathfrak{a})$ again denote a solution to (1.2). This subsection says more about the curvature of A in small radius balls in X the centers and radii of the balls are described by the first paragraph in Section 5c. (This subsection is not used in the proof of Theorem 1.2.)

It is assumed at the outset that $\mu$ is chosen from $(0, \frac{1}{100}]$ and that $c$ is chosen large enough to invoke Propositions 5.1, 5.2 and now add the extra condition that $c$ be large enough to invoke Proposition 5.4 also. Let p denote again the center of the ball in question and let r denote its radius. By assumption, $r \le r_{c_\mu}$. It is also assumed that $r \le \frac{1}{c}$ but that $rrK(r)$ is greater than $c_\mu$ so that Propositions 5.2 and 5.4 can be called on.

The self dual part of A's curvature, $F_A^+$, obeys the pointwise bound

$$|F_A^+| \le c_\mu \frac{1}{c} e^{-\sqrt{rrK(r)}/\kappa_\mu} \frac{1}{r^2}$$

(5.37)

on $B_{(1-\mu)r}$. This is because of (5.3) (and the $x^2 e^{-x} \le 16 e^{-x/2}$ bound with $x = \sqrt{rrK(r)}/c_\mu$.) Meanwhile, Proposition 5.5 bounds the $L^2$ norm on $B_{(1-3\mu)r}$ of the part of $F_A^-$ that is pointwise orthogonal to $\sigma_\ddagger$. The part of $F_A^-$ along $\sigma_\ddagger$ is not well constrained however. The focus can be put on this part of $F_A^-$ by introducing a new connection on the bundle P over $B_{(1-\mu)r}$ that is denoted by $\hat{A}$ and defined by the rule:

$$\hat{A} = A - \tfrac{1}{4}[\sigma_\ddagger, d_A\sigma_\ddagger] .$$

(5.38)

This connection $\hat{A}$ is defined in part so that $\nabla_{\hat{A}}\sigma_\ddagger = 0$. This is equivalent to the following: An isomorphism over $B_{(1-\mu)r}$ between the product principal SU(2) (or SO(3)) bundle and P that identifies $\sigma_\ddagger$ with a constant matrix (to be denoted by $\tau$) writes $\hat{A}$ as $\theta_0 + A\tau$ with $\theta_0$ denoting the product connection and with A denoting an $\mathbb{R}$-valued 1-form.

The curvature of $\hat{A}$ is proportional to $\sigma_\ddagger$ and it obeys

$$\langle \sigma_\ddagger F_{\hat{A}} \rangle = \langle \sigma_\ddagger F_A \rangle - \tfrac{1}{4}\langle \sigma_\ddagger d_A\sigma_\ddagger \wedge d_A\sigma_\ddagger \rangle .$$

(5.39)

This is the formula for the 2-form dA. It follows as a consequence of this identity (and the second bullet of Proposition 5.2 and (5.28) and Proposition 5.4) that

$$\int_{B_{(1-3\mu)r}} |F_{\hat{A}} - F_A|^2 \le c_\mu \frac{1}{c^2} e^{-\sqrt{rrK(r)}/c_\mu} (1 + \frac{F^2(r)}{r^2 K^2(r)}) .$$

(5.40)

What with (5.37), this implies that $\langle \sigma_\ddagger F_{\hat{A}} \rangle$ is nearly self-dual on $B_{(1-3\mu)r}$ (assuming that $\frac{F^2(r)}{r^2 K^2(r)}$ is not too big and $rrK(r)$ is not too small:

$$\int_{B_{(1-3\mu)r}} |F_{\hat{A}}^+|^2 \le c_\mu \frac{1}{c^2} e^{-\sqrt{rrK(r)}/c_\mu} (1 + \frac{F^2(r)}{r^2 K^2(r)})$$

(5.41)

This bound implies in turn that $\hat{A}$ can be modifed by the addition of a relatively small term (proportional to $\sigma_\ddagger$) so that the resulting connection has anti-self dual curvature. The next paragraphs explains how to do this.

Standard elliptic theory finds a section of $\Lambda^+$ over $B_{(1-3\mu)r}$ (to be denoted by $u$) that obeys the conditions

- $d^\dagger du = \langle \sigma_\ddagger F_{\hat{A}}{}^+ \rangle$
- $u = 0$ on $\partial B_{(1-3\mu)r}$.

(5.42)

Give these equations, then the bounds given below in (5.43) are straightforward consequences of the bounds in (5.41) and (1.9), and the fact that $d^\dagger d = \nabla^\dagger \nabla + \mathfrak{R}$ with $\mathfrak{R}$ denoting the Riemann curvature endomorphism that appears in (2.1).

- $\int_{B_{(1-4\mu)r}} |\nabla u|^2 \le c_\mu \frac{1}{c^2} e^{-\sqrt{rrK(r)}/c_\mu} (1 + \frac{F^2(r)}{r^2 K^2(r)}) r^2$
- $\int_{B_{(1-4\mu)r}} |\nabla \nabla u|^2 \le c_\mu \frac{1}{c^2} e^{-\sqrt{rrK(r)}/c_\mu} (1 + \frac{F^2(r)}{r^2 K^2(r)})$ .

(5.43)

Let $\mathcal{A}$ denote the connection $\hat{A} + *du\sigma_\ddagger$. Then $\nabla_\mathcal{A} \sigma_\ddagger = 0$ so $\sigma_\ddagger$ is $\mathcal{A}$-covariantly constant. This implies that $F_\mathcal{A}$ is proportional to $\sigma_\ddagger$ and that $\langle \sigma_\ddagger F_\mathcal{A} \rangle$ is a closed 2-form. Meanwhile, by virtue of (5.42), the curvature of $\mathcal{A}$ is anti-self dual (which is to say that $F_\mathcal{A}{}^+ = 0$). And, by virtue of (5.43):

$$\int_{B_{(1-4\mu)r}} |F_\mathcal{A} - F_A|^2 \le c_\mu \frac{1}{c^2} e^{-\sqrt{rrK(r)}/c_\mu} (1 + \frac{F^2(r)}{r^2 K^2(r)}) .$$

(5.44)

Thus, supposing that $\frac{F^2(r)}{r^2 K^2(r)}$ is not too large and $rrK(r)$ is not too small, the connection A is almost reducible and it almost has anti-self dual curvature.

As a parenthetical remark, the fact that $\langle \sigma_\ddagger F_\mathcal{A} \rangle$ is closed and anti-self dual implies the following (see Step 6 in Section 5c of [T3]): Suppose that $q \in B_{(1-4\mu)r}$ and that $\rho$ is positive but less than $\text{dist}(q, \partial B_{(1-4\mu)r})$. Then

$$\int_{\text{dist}(\cdot, q) \le \rho} |F_\mathcal{A}|^2 \le (\tfrac{\rho}{r})^4 (1 - c_0 r^2) \int_{B_r} |F_\mathcal{A}|^2 .$$

(5.45)

For what it is worth, the $q = p$ version of this inequality and the version of (5.44) with $r = \rho$ lead to the following observation: If $r \ge c_\mu \frac{1}{\rho K(\rho)} (\ln(\rho K(\rho))^2$, then

$$F(\rho) \le (\tfrac{\rho}{r})^4 (1 - c_0 r^2) F(r) + c_\mu \frac{1}{c^2} e^{-\sqrt{r\rho K(\rho)}/c_\mu} .$$

(5.46)

This is a uniform bound on F(ρ) supposing the existence of one for F(r)

## 6. Proof of Theorem 1.2

The preceding sections supplied the tools that are used in the upcoming Sections 6a-6e to prove the assertions of Theorem 1.2. By way of a reminder, the input for Theorem 1.2 consists of a sequence $\{r_n, (A_n, a_n)\}_{n \in \mathbb{N}}$ with $\{r_n\}_{n \in \mathbb{N}}$ denoting an increasing, unbounded sequence of positive real numbers, and with each $n \in \mathbb{N}$ version of $(A_n, a_n)$ obeying the $r = r_n$ version of (1.2). This input data is taken for granted in what follows.

### a) $L^2_1$ and $L^\infty$ limits

The proof of Theorem 1.1 starts here with the extraction of a subsequence from $\{|a_n|\}_{n \in \mathbb{N}}$ that has certain desirable properties, one being its weak convergence in the $L^2_{1;loc}$ topology on X to an $L^\infty$ function. This subsequence is described in Proposition 6.1. This proposition is the Vafa-Witten analog of Proposition 2.1 in [T2] and Proposition 10.1 in [T3].

**Proposition 6.1**: *There exists $\kappa > 1$ with the following significance: Let $\{r_n\}_{n \in \mathbb{N}} \subset (1, \infty)$ denote an increasing and unbounded sequence; and for each $n \in \mathbb{N}$, let $(A_n, a_n)$ denote a solution to the $r = r_n$ version of (1.2). There exists a subsequence $\Lambda \subset \mathbb{N}$ such that the bulleted items listed below hold. This list uses $\underline{X}$ to denote X when X is compact, and to denote any given open set in X with compact closure and smooth boundary otherwise.*

- *The sequences $\{ \int_{\underline{X}} (|d|a_n||^2 + |a_n|^2) \}_{n \in \Lambda}$ and $\{\sup_{\underline{X}} |a_n|\}_{n \in \Lambda}$ are bounded.*

- *The sequence $\{|a_n|\}_{n \in \Lambda}$ converges weakly in the $L^2_1$ topology on $\underline{X}$ and strongly in all $p < \infty$ versions of the $L^p$ topology on $\underline{X}$.*

- *The limit function (it is denoted by $|v_\ddagger|$) is an $L^\infty$ function whose value can be defined at each point in $\underline{X}$ by the rule whereby $|v_\ddagger|(p) = \limsup_{n \in \Lambda} |a_n|(p)$ for each $p \in \underline{X}$.*

- *The sequence $\{\langle a_n \otimes a_n \rangle\}_{n \in \Lambda}$ converges strongly in any $q < \infty$ version of the $L^q$ topology on the space of sections of $\Lambda^+ \otimes \Lambda^+$ over $\underline{X}$. The limit section is denoted by $v_\ddagger \otimes v_\ddagger$ and its trace is the function $|v_\ddagger|^2$.*

- *Use f to denote a given $C^0$ functional The sequences $\{ \int_{\underline{X}} f |\nabla_{A_n} a_n|^2 \}_{n \in \Lambda}$ and $\{r_n^2 \int_{\underline{X}} f |[a_n; a_n]|^2 \}_{n \in \Lambda}$ converge. The limit of the first sequence is denoted by $Q_{\nabla, f}$ and that of the second by $Q_{[\cdot], f}$. These are such that*

$$\tfrac{1}{2}\int_X d^*df\,|v_\ddagger|^2 + Q_{\nabla,f} + Q_{[\cdot],f} + \int_X f\langle\mathfrak{R}, v_\ddagger\otimes v_\ddagger\rangle = 0,$$

- *Fix $p \in \underline{X}$ and let $G_p$ denote the (Dirichelet) Green's function on $\underline{X}$ with pole at p for the operator $d^\dagger d + 1$. The sequence that is indexed by $\Lambda$ with n'th term being the integral of $G_p(|\nabla_{A_n}a_n|^2 + r_n^2|[a_n; a_n]|^2)$ is bounded. Let $Q_{\Diamond,p}$ denote the lim-inf of this sequence of integrals. The function $|v|^2$ obeys the equation*

$$\tfrac{1}{2}|v_\ddagger|^2(p) + Q_{\Diamond,p} = \int_{\underline{X}} G_p\left(\tfrac{1}{2}|v_\ddagger|^2 - \langle\mathfrak{R}, v_\ddagger\otimes v_\ddagger\rangle\right).$$

The proof of this proposition differs only cosmetically from Part II of the proof of Proposition 2.1 in [T] to which the reader is directed. The proof is not given here.

**b) Hölder continuity across the zero locus**

The third bullet of the Proposition 10.1 defines the values of $|v_\ddagger|$ at each point in X. Granted this definition, let Z denote the set of points in X (or $\underline{X}$ when X is not compact) where $|v_\ddagger| = 0$. Thus, p being in Z means that

$$\limsup_{n\in\Lambda} |a_n|(p) = 0.$$

(6.1)

The next proposition makes the formal assertion to the effect that Z is a closed set in X and that $|v_\ddagger|$ is uniformly Hölder continuous across Z.

**Proposition 6.2**: *There exists $\kappa > 1$ with the following significance: Let $\{r_n, (A_n, a_n)\}_{n\in\mathbb{N}}$ denote the input sequence for Proposition 6.1, and let $\Lambda$ denote the subsequence of $\mathbb{N}$ given in Proposition 6.1. Use the corresponding sequence $\{a_n\}_{n\in\Lambda}$ to define the function $|v_\ddagger|$ as per the instructions in Proposition 6.1. Define Z to be the set of points where the function $|v_\ddagger|$ is zero, which is where (6.1) holds. This set Z is a closed set in X. And, if $p \in Z$, then $|v_\ddagger|$ on the radius $\kappa^{-1}$ ball centered at p obeys $|v_\ddagger| \leq \kappa\,\mathrm{dist}(p,\cdot)^{1/\kappa}$.*

This proposition is proved momentarily given a second proposition that concerns an (approximate) analog of the function $\mathcal{K}$ that is defined using $|v_\ddagger|$ in lieu of $|a|$. This function is defined after choosing a point $p \in X$ (in $\underline{X}$ if X is not compact). Given p, the corresponding function is defined on $(0, r_0]$. It is denoted by $\mathcal{K}$ and it is given by the rule

$$r \to \mathcal{K}(r)^2 = \tfrac{1}{r^3}\int_{\partial B_r} |v_\ddagger|^2.$$

(6.2)

The function $\mathcal{K}$ is bounded because $|v_\ddagger|$ is bounded. (Proposition 6.1's third bullet says in part that $|v_\ddagger|$ is bounded.) Note that there is a precise analog of $\mathcal{K}$ in [T3], defined in (10.2) of this reference. Section 10a of [T3] also has a precise analog of the following proposition about $\mathcal{K}$.

**Proposition 6.3**: *There exists $\kappa > 1$ with the following significance: Let $\{r_n, (A_n, a_n)\}_{n \in \mathbb{N}}$ denote the input sequence for Proposition 6.1, and let $\Lambda$ denote the subsequence of $\mathbb{N}$ given in Proposition 6.1. Use the the instructions in Proposition 6.1 to define the function $|v_\ddagger|$ from the sequence $\{a_n\}_{n \in \Lambda}$. Define $Z$ to be the zero locus of $|v_\ddagger|$. Supposing that $p \in Z$ (in $Z \cap \underline{X}$ if $X$ is not compact), and that $r \in (0, r_0]$, then $\mathcal{K}(r) \leq \kappa r^{1/\kappa}$.*

 This proposition is proved momentarily. Assume it to be true, then the argument in Section 10a of [T3] that proves the latter's Proposition 10.2 given the latter's Proposition 10.3 can be borrowed virtually verbatim to prove Proposition 6.2 here from Proposition 6.3. Since this borrowed argument is short, it is reproduced below.

*Proof of Proposition 6.2*: The assertion that $Z$ is closed follows from the Hölder bound for the function $|v_\ddagger|$ because the Hölder bound implies that if ever $\{p_i\}_{i \in \mathbb{N}}$ is a sequence in $X$ that converges to a point in $Z$, then $\lim_{i \to \infty} |v_\ddagger|(p_i) = 0$.
 To prove the Hölder bound, fix $p \in Z$, and then $r \in (0, r_0]$. Having done this, use the canonical bump function $\chi$ to construct a function on $X$ to be denoted by $\chi_\diamond$ that is equal to 0 where the distance to $p$ is greater than $\frac{27}{32} r$ and it is equal to 0 where the distance to $\mu$ is less than $\frac{26}{32} r$. This function can and should be consructed so that $|d\chi_\diamond| \leq c_0 \frac{1}{r}$ and $|\nabla d\chi_\diamond| \leq c_0 \frac{1}{r^2}$. Supposing that $q \in B_{3r/4}$, let $G_q$ again denote the Dirichelet Greens function for the operator $d^*d$ on $B_r$ with pole at the point $q$. (See (1.10).) Use $f$ to be the function $(1 - \chi_\diamond) G_q$. Use this function in the fifth bullet of Proposition 10.1 and use $q$ in lieu of $p$ in the sixth bullet of Proposition 10.1. Subtract the fifth bullet's identity from sixth bullet's identity to obtain an inequality for $|v|^2(q)$ that takes the form

$$\tfrac{1}{2}|v_\ddagger|^2(q) \leq \int_X \chi_{1/4} G_q (\tfrac{1}{2}|v_\ddagger|^2 - \langle \mathfrak{R}, v_\ddagger \otimes v_\ddagger \rangle) + \int_X (\tfrac{1}{2} d^* d\chi_{1/4} G_p - \langle d\chi_{1/4}, dG_q \rangle)|v_\ddagger|^2 .$$

(6.3)

This in turn leads to the bound

$$|v_\ddagger|^2(q) \leq c_0 r^2 + c_0 \frac{1}{r^4} \int_{B_r} |v_\ddagger|^2 .$$

(6.4)

To explain: The $c_0 r^2$ term in (6.4) accounts for the left most integral on the right hand side of (6.3) because $|v_{\ddagger}| \leq c_0$ and because of the middle bullet in (1.9). The right most integral on the right hand side of (6.3) is accounted for by the product of $c_0 \frac{1}{r^4}$ times the integral of $|v_{\ddagger}|^2$ that appears on the right hand side of (6.4). Proposition 6.3 says that the right hand side of (6.4) is bounded by $c_0 r^{1/c_0}$ when $r \leq c_0^{-1}$. This gives the Hölder norm bound of Proposition 6.2.

***Proof of Proposition 6.3***: Fix $\mu \in (0, \frac{1}{100}]$ and then fix $c$ sufficiently large (the lower bound is $c_\mu$) so that the Proposition 4.1 and 5.1 can be invoked. Suppose that $p \in Z$ (and in $\underline{X}$ if X is not compact). (This is to say that (6.1) holds). For each $n \in \Lambda$, the point p has the version of $r_{c\mu}$ that is defined by the pair $(A_n, a_n)$. Denote this version of $r_{c\mu}$ by $r_{c\mu,n}$. There are two instances to consider with regards to the sequence $\{r_{c\mu,n}\}_{n \in \Lambda}$. In the first instance, $\{r_{c\mu,n}\}_{n \in \Lambda}$ has a subsequence that is bounded away from zero; and in the second instance, there is no such subsequence. These are treated in Parts 1 and 2 of the proof.

*Part 1*: Suppose that $\{r_{c\mu,n}\}_{n \in \Lambda}$ has a subsequence that is bounded away from zero. Let $\Theta$ denote this subsequence and let $r_\diamond$ denote $\frac{1}{2} \lim_{n \to \infty} r_{c\mu,n}$. If n is sufficiently large, the number $r_\diamond$ will be less than $(1 - 6\mu) r_{c\mu,n}$. Let $\kappa_n$ denote the version of the function $\kappa$ that is defined by the pair $(A_n, a_n)$. If $n \in \Theta$ is large, and if $c$ is greater than $c_\mu$, then the first bullet of Proposition 5.1 implies that $|a_n| \geq \frac{1}{2\pi} \kappa_n(r_\diamond)$ in the radius $r_\diamond$ ball centered at p. What with (6.1), this is possible only in the event that $\lim_{n \in \Theta} \kappa_n(r_\diamond) = 0$. But, if this is so, then $\kappa(r) = 0$ for any $r < r_\diamond$ (a consequence of the second bullet in Proposition 6.1 and the fact that each $\kappa_n$ is a non-decreasing function). Since the constant function 0 is less than $cr^{1/c}$ for an $c > 1$, Proposition 6.3's assertion holds in this case.

*Part 2*: Suppose now that $\lim_{n \to \infty} r_{c\mu,n} = 0$. Proposition 4.1 introduces the number $\lambda_{c\mu}$. Fix $r_\diamond$ to be less than $\lambda_{c\mu}^{-1} r_0$. For each $n \in \Lambda$, let $\kappa_n$ again denote the version of the function $\kappa$ that is defined by the pair $(A_n, a_n)$. There are now two subcases to distinguish. In the first case, there exists a subsequence of $\Lambda$ to be denoted by $\Theta$ where $\lim_{n \in \Theta} r_n r_\diamond \kappa_n(r_\diamond) < 1$. Since $r_\diamond$ is fixed, this occurs only in the event that $\lim_{n \to \infty} \kappa_n(r_\diamond) \to 0$. The argument in Part 1 can be repeated to conclude that $\kappa(r) = 0$ for all $r \leq r_\diamond$ which is consistent with the assertion of Proposition 6.3.

In the second case, for any choice of $r_\diamond$, there exists a positive integer $n_\diamond$ such that $r_n r_\diamond \kappa_n(r_\diamond) \geq 1$ when n is greater than $n_\diamond$. By taking $n_\diamond$ larger if necessary, $r_{c\mu,n}$ will be smaller than $\frac{1}{100} r_\diamond$ when $n > n_\diamond$. Take $n_\diamond$ so that this last bound is true also. With the preceding understood, suppose that $n > n_\diamond$. The $(A_n, a_n)$ version of Proposition 4.1 says that $K_n(r_2) \leq (\frac{r_2}{r_1})^{1/\lambda_{c_\mu}} K_n(r_1)$ if $r_1 \in (r_\diamond, \lambda_{c_\mu}^{-1} r_0]$ and $r_2 \in [r_\diamond, r_1]$.

By virtue of the second bullet in Proposition 6.1, the bound just stated implies the following with regards to the function $\mathcal{K}$ (since $r_\diamond$ can be as small as desired): For any given $r_1 \in (0, \lambda_{c_\mu}^{-1} r_0]$ and $r_2 \in (0, r_1]$, the function $\mathcal{K}$ obeys

$$\mathcal{K}(r_2) \leq (1 + c_0 r_1^2) \left(\frac{r_2}{r_1}\right)^{1/\lambda_{c_\mu}} \mathcal{K}(r_1).$$

(6.5)

The $r_2 = r$ and $r_1 = (\frac{r_2}{r_1})^{1/\lambda_{c_\mu}} r_0$ version of this gives Proposition 6.3's asserted bound.

**b) The sequence $\{r_{c_\mu,n}\}$ where $|v_\ddagger| > 0$**

Fix $\mu \in (0, \frac{1}{100}]$ and then fix $c > c_\mu$ (chosen so that the propositions in Sections 4 and 5 can be invoked.) Let p denote a chosen point in X (in $\underline{X}$ if X is not compact) with $|v_\ddagger|(p) > 0$. The lemma below makes a formal assertion to the effect that the sequence $\{r_{c_\mu,n}\}_{n \in \Lambda}$ has strictly positive lim-inf.

**Lemma 6.4**: *There exists $\kappa > 1$, and given $\mu \in (0, \frac{1}{100}]$, there exists $\kappa_\mu > \kappa$ that depends only on $\mu$; and these numbers have the following significance: Let $\{r_n, (A_n, a_n)\}_{n \in \Lambda}$ and $|v_\ddagger|$ denote the output data from Proposition 6.1. Fix $c > \kappa_\mu$ and a point $p \in X$ (in $\underline{X}$ if X is not compact) where $|v_\ddagger|(p) > 0$. For each $n \in \Lambda$, reintroduce the number $r_{c_\mu,n}$ for p (which is the version of $r_{c_\mu}$ that is defined by $(A_n, a_n)$). Then $\liminf_{n \in \Lambda} r_{c_\mu,n} > 0$.*

*Proof of Lemma 6.4*: Suppose to the contrary that there is a subsequence $\Xi \subset \Lambda$ with the corresponding sequence $\{r_{c_\mu,n}\}_{n \in \Xi}$ having limit zero. The three steps that follow generate nonsense from this observation.

Step 1: Let $\Theta_p \subset \Lambda$ denote a subsequence with $\lim_{n \in \Theta_p} |a_n|(p) = |v_\ddagger|(p)$. Given $\varepsilon > 0$, there is an integer $n_\varepsilon$ such that $|a_n|(p) \geq (1 - \varepsilon)|v_\ddagger|(p)$ when n from $\Theta$ is greater than $n_\varepsilon$. This then implies that $K_n(r) \geq \sqrt{2\pi}(1 - \varepsilon)|v_\ddagger|(p)$ for each $r \in (0, r_0]$ when n from $\Theta$ is greater than $n_\varepsilon$. This implies that $\mathcal{K}(r) \geq \sqrt{2\pi}(1 - \varepsilon - c_0 r^2)|v_\ddagger|(p)$ for each $r \in (0, r_0]$ (by virtue of the second bullet of Proposition 6.1.).

Step 2: Fix $r \in (0, r_0]$ and $\varepsilon > 0$. If n is from $\Xi$ and n is large (the lower bound is determined by r and $\varepsilon$), then $K_n(r) \geq (1 - \varepsilon) \mathcal{K}(r)$. This is due to the second bullet in Proposition 6.1 also). Therefore, when n is large, $r_n r K_n(r) > 1$. As a consequence, what is said by Proposition 4.1 can be applied using this value of r and $(A_n, a_n)$. In particular, if n is so large that $r_{c_\mu,n} < r$, then the $r_2 = r$ and $r_1 = \lambda_{c_\mu}^{-1} r_0$ version of the proposition's inequality says that $K_n(r) \leq (1+c_0 r_1^2) (\frac{r}{r_1})^{1/\lambda_{c_\mu}} K_n(r_1)$.

<u>Step 3</u>: Since $K_n(r_1) \leq c_0$ (by virtue of Lemma 2.1), and since $K_n(r) \geq (1-\varepsilon) \mathcal{K}(r)$, and since $\mathcal{K}(r) \geq c_0^{-1} |v_\ddagger|(p)$, and since $|v_\ddagger|(p) > 0$, the conclusions in Step 2 lead to nonsense if $r_{c_\mu,n}$ is less than $(c_\mu |v|(p))^{\lambda_{c_\mu}} \lambda_{c_\mu}^{-1}$.

### c) Pointwise and $C^0$ convergence

The next proposition makes a formal assertion to the effect that the sequence $\{|a_n|\}_{n\in\Lambda}$ from Proposition 6.1 converges pointwise to $|v_\ddagger|$ in the $C^0$ topology (on $\underline{X}$ if X is not compact).

**Proposition 6.5**: *Let $\{r_n, (A_n, a_n)\}_{n\in\Lambda}$ and $|v_\ddagger|$ denote the data supplied by Proposition 6.1. The function $|v_\ddagger|$ is continuous and the sequence $\{|a_n|\}_{n\in\Lambda}$ converges in the $C^0$ topology to $|v_\ddagger|$ on compact subsets of X.*

*Proof of Proposition 6.5*: The proof has two parts. The first part proves pointwise convergence of $\{|a_n|\}_{n\in\Lambda}$ to $|v_\ddagger|$, and the second part proves that $|v_\ddagger|$ is continuous and that $\{|a_n|\}_{n\in\Lambda}$ converges uniformly to $|v_\ddagger|$ in the $C^0$ topology on compact subsets of X.

*Part 1*: Fix $p \in X$ (in $\underline{X}$ if X is not compact.) This step proves that the sequence of numbers $\{|a_n|(p)\}$ converges to $|v_\ddagger|(p)$. If $|v_\ddagger|(p) = 0$, this follows by virtue of the definition of $|v_\ddagger|(p)$ in the third bullet of Proposition 6.1. With this understood, assume henceforth that $|v_\ddagger|(p) \neq 0$.

Fix $\mu \in (0, \frac{1}{100}]$ and then $c > c_\mu$ with the lower bound chosen so that the conclusions of Propositions 5.1 and 5.2 and Lemma 6.4 can be invoked Given a positive integer k, invoke Lemma 6.4 with c replaced by ck. For each $n \in \Lambda$, let $r_{kc_\mu,n}$ denote the kc version of what Lemma 6.4 denotes by $r_{c_\mu,n}$. According to Lemma 6.4, the sequence $\{r_{kc_\mu,n}\}_{n\in\Lambda}$ has a positive lower bound. Use $r_{0k}$ to denote half of such a lower bound. Fix any positive r that is less than $\frac{4}{3} r_{0k}$ and less than $\frac{1}{kc}$. If $n \in \Lambda$, then the $a = a_n$ version of the top bullet in Proposition 5.1 says that

$$||a_n| - \tfrac{1}{\sqrt{2\pi}} K_n(r)| \le c_\mu \tfrac{1}{kc}$$

(6.6)

in the radius $(1-\mu)r$ ball centered at p.

Reintroduce the function $\mathcal{K}$ on $(0, r_0)$ from (6.2). Since $\lim_{n\to\infty} K_n(r) = \mathcal{K}(r)$ for any positive r (by virtue of Proposition 6.1), the inequality in (6.6) implies in turn that

$$||a_n| - \tfrac{1}{\sqrt{2\pi}} \mathcal{K}(r)| \le c_\mu \tfrac{1}{kc}$$

(6.7)

in the radius $(1-\mu)r$ ball centered at p if n is sufficiently large. (If r is the minimum of $r_{0k}$ and $\tfrac{1}{2kc}$, then (6.7) holds for n greater than a number that depends only on k.)

To see that (6.7) implies pointwise convergence, fix a subsequence $\Theta_p \subset \Lambda$ such that $\lim_{n\in\Theta_p} |a_n|(p) = |v_\ddagger|(p)$. Use integers from $\Theta_p$ in (6.7) and evaluate $|a_n|$ at p to see that

$$||v_\ddagger|(p) - \tfrac{1}{\sqrt{2\pi}} \mathcal{K}(r)| \le c_\mu \tfrac{1}{kc}$$

(6.8)

when r is the minimum of $r_{0k}$ and $\tfrac{1}{2kc}$. Now use this bound with (6.7) (with no constraint on n but for it being large) to see that

$$||a_n| - \tfrac{1}{\sqrt{2\pi}} |v_\ddagger|(p)| \le c_\mu \tfrac{1}{kc}$$

(6.9)

on this same ball of radius r ball centered at p if n ∈ Λ is sufficiently large given k. Since k can be any positive integer, (6.9) at the point p proves that $\lim_{n\in\Lambda} |a_n|(p) = |v_\ddagger|(p)$.

*Part 2*: The fact that $|v_\ddagger|$ is continuous is also a consequence of (6.9). To see this, take a point q in the radius r ball centered at p (with r still the minimum of $r_{0k}$ and $\tfrac{1}{2kc}$). Fix a sequence $\Theta_q \subset \Lambda$ so that $\{|a_n|\}_{n\in\Theta_q}$ converges to $|v_\ddagger|(q)$. Then (6.9) implies that

$$||v_\ddagger|(q) - |v_\ddagger|(p)| \le c_\mu \tfrac{1}{kc} .$$

(6.10)

Since this hold for any point q in the radius r ball centered at p; and since k can be any positive integer, this inequality proves that $|v_\ddagger|$ is continuous at points in X−Z (because p is assumed to be in X−Z). Proposition 6.2 proves that $|v_\ddagger|$ is also continuous across Z.

To consider the question of uniform $C^0$ convergence, fix a positive integer, k, and note that each point in X has its own version of $r_{0k}$. Therefore, supposing that U ⊂ X−Z is a compact set, there is a finite set of balls that cover U with the radius of each ball being the minimum of $r_{\Diamond k}$ at its center point and $\frac{1}{2kc}$. Fix such a cover. With the cover fixed, and given positive $\varepsilon$ (and k), there exists a number $n_{\varepsilon,k}$ with the following property: Let p denote the center point of a ball in the cover and let $r_p$ denote the minimum of p's version of $r_{\Diamond k}$ and $\frac{1}{2kc}$. If $n \geq n_{\varepsilon,k}$, then $|\mathcal{K}_n(r_p) - \mathcal{K}(r_p)| \leq \varepsilon$. Indeed, the existence of $n_{\varepsilon,k}$ follows from the second bullet of Proposition 6.1. This bound and the inequalities (6.6)-(6.10) as applied to each of the finite balls in U's cover imply that

$$||a_n| - |\nu_{\ddagger}|| \leq c_\mu (\varepsilon + \frac{1}{kc})$$

(6.11)

on U if $n \geq n_{\varepsilon,k}$. This last inequality implies the uniform convergence assertion of Proposition 6.5 because if positive $\varepsilon$ is chosen a priori, then (6.11) can be invoked with k being the smallest integer greater than $\varepsilon^{-1}c$.

### d) The self-dual 2-form $\nu_{\ddagger}$ on small radius balls where $|\nu_{\ddagger}| > 0$

The proposition that follows constructs Theorem 1.2's self-dual 2-form $\nu_{\ddagger}$ on small radius balls around the points in X where $|\nu_{\ddagger}|(p) > 0$.

**Proposition 6.6**: *There exists $\kappa > 1$, and given $\mu \in (0, \frac{1}{100}]$, there exists $\kappa_\mu > \kappa$ that depends only on $\mu$; and these numbers have the following significance: Let $\{r_n, (A_n, a_n)\}_{n \in \Lambda}$ and $|\nu_{\ddagger}|$ denote the output data from Proposition 6.1. Fix $c > \kappa_\mu$ and a point $p \in X$ (in $\underline{X}$ if X is not compact) where $|\nu_{\ddagger}|(p) > 0$. Define the number $r_{\Diamond}$ to be half of the minimum of $\liminf_{n \in \Lambda} r_{c_\mu,n}$ and $\frac{1}{c}$. For each $n \in \Lambda$, the conclusions of Propositions 5.1 and 5.2 hold with $(A, a) = (A_n, a_n)$ and $r = r_{\Diamond}$ and the given value of c. In particular, the section $a_n$ can be written over $B_{r_{\Diamond}}$ as $a_n = \nu_n \sigma_{\ddagger n} + \mathfrak{a}_n$ as explained in Proposition 5.2. After multiplying some elements of $\{\nu_n\}_{n \in \Lambda}$ by -1 (perhaps), the resulting version of the sequence $\{\nu_n\}_{n \in \Lambda}$ on $B_{r_{\Diamond}}$ converges strongly in the $L^2_1$ and $C^0$ topologies to a smooth, closed, self-dual 2-form on $B_{r_{\Diamond}}$. Moreover, $\{\nu_n\}_{n \in \Lambda}$ converges in the $L^2_2$ topology to this same 2-form if there is an a priori upper bound for the sequence whose n'th member is the $L^2$ norm of $\frac{1}{r_n} F_{A_n}$ on the ball $B_{(1+\mu)r_{\Diamond}}$.*

With regards to the multiplications of elements in $\{\nu_n\}_{n \in \Lambda}$ by -1: The depiction of a in Proposition 5.2 as $a = \nu \sigma_{\ddagger} + \mathfrak{a}$ defines the section $\nu$ (a section of $\Lambda^+$) only up to multiplication by -1. This is because there is an automorphism of the bundle P over Proposition 5.2's ball $B_r$ that changes $\sigma_{\ddagger}$ and $-\sigma_{\ddagger}$. This $\pm 1$ ambiguity is why some elements in the sequence $\{v_n\}_{n\in\Lambda}$ from Lemma 6.4 are multiplied by -1 before the strong convergence assertion for the whole sequence can be made.

***Proof of Proposition 6.6***: The proof of the proposition has four parts.

*Part 1*: It follows from the $r = r_\diamond$ version of Proposition 5.2 that the sequence $\{v_n\}_{n\in\Lambda}$ has bounded $L^2_1$ norm and bounded $L^\infty$ norm on the radius $\frac{9}{8} r_\diamond$ ball centered at p. This implies that there is a subsequence that converges weakly in the $L^2_1$ topology to a limit which is a priori an $L^2_1$ and $L^\infty$ section of $\Lambda^+$ over this same ball. By virtue of the Rellich lemma, this convergence is strong in the $L^2$ topology. However, by virtue of the fourth bullet of Proposition 6.1 and the fourth bullet of Proposition 5.2 (applied to $(A_n, a_n)$ with r being $\frac{9}{8} r_\diamond$), the sequence $\{v_n\}_{n\in\Lambda}$ can be changed by multiplying each member by either 1 or -1 so that the whole of the resulting sequence converges strongly in the $L^2$ topology on the radius $\frac{9}{8} r_\diamond$ ball centered at p. Make these +1 or -1 multiplications and now call the result $\{v_n\}_{n\in\Lambda}$. The strong $L^2$ convergence of the whole sequence implies in turn that the whole sequence now converges weakly in $L^2_1$ (as opposed to just a subsequence converging) on this same ball.

*Part 2*: Let $v_\ddagger$ denote the limit of the sequence $\{v_n\}_{n\in\Lambda}$. It is a consequence of Bullets 2-4 of Proposition 5.2 that $dv_\ddagger = 0$. This implies that $v_\ddagger$ is harmonic and therefore it is a smooth section of $\Lambda^+$. Moreover, since the norm of $v_\ddagger$ is Proposition 6.1's function $|v_\ddagger|$, it follows that the latter function is smooth on X−Z.

By virtue of those same bullets in Proposition 5.2, the sequence $\{d(v_n - v_\ddagger)\}_{n\in\Lambda}$ converges to zero in the $L^2$ topology on the radius $\frac{9}{8} r_\diamond$ ball centered at p. This fact is used directly to prove that $\{v_n\}_{n\in\Lambda}$ converges to $v_\ddagger$ strongly in the $L^2_1$ topology on $B_{r_\diamond}$. To prove this, use the function $\chi$ to construct a smooth function on $B_{7r_\diamond/8}$ that is equal to 1 on $B_{r_\diamond}$ and has compact support. This function can and should be constructed so that the norm of its derivative is bounded by $c_0 \frac{1}{r_\diamond}$ and so that the norm of its Hessian is no greater than $c_0 \frac{1}{r_\diamond^2}$. Denote this function by $\chi_\diamond$. A back-and-forth sequence of integration by parts can be used to see that

$$\int_{B_{17r_\diamond/16}} |\nabla(v_n - v_\ddagger)|^2 \leq c_0 \int_{B_{9r_\diamond/8}} \chi_\diamond^2 |d(v_n - v_\ddagger)|^2 + c_0 \frac{1}{r_\diamond^2} \int_{B_{9r_\diamond/8}} |v_n - v_\ddagger|^2 .$$

(6.12)

Since the $n \to \infty$ limit of the right hand side of (6.6) is zero (with n from $\Lambda$), this is necessarily the case for the $n \to \infty$ limit of the left hand side.

*Part 3*: To prove the $C^0$ convergence of $\{v_n\}_{n \in \Lambda}$ to $v_\ddagger$, first fix $n \in \Lambda$ to consider the $(A_n, a_n)$ version of (5.15) on the radius $\frac{17}{16} r_\diamond$ ball centered at p. If $c \geq c_\mu$, then the equations in (5.15) (and the fact that $dv = 0$) lead to an equation on this ball of the form

$$\nabla^\dagger \nabla (v_n - v_\ddagger) = \mathfrak{e}_n \quad (6.13)$$

with $\mathfrak{e}_n$ being a small term. To say that $\mathfrak{e}_n$ is small means in particular that

$$\lim_{n \in \Lambda} \sup_{q \in B_{33r_\diamond/32}} \int_{B_{17r_\diamond/16}} \frac{1}{\text{dist}(q, \cdot)^2} |\mathfrak{e}_n|^2 = 0 \quad (6.14)$$

This limit follows from Proposition 5.2's third and fourth bullets and from the $L^2_1$ convergence of $\{v_n\}_{n \in \Lambda}$ to $v_\ddagger$.

Given a point q with distance less than $r_\diamond$ from p, let $\mathcal{G}_q$ denote the Dirichlet Green's function on $B_{r_\diamond}$ with pole at q for the operator $\nabla^\dagger \nabla$ acting on sections of $\Lambda^+$. The norm of this Green's function is at most $c_0 \frac{1}{\text{dist}(q, \cdot)^2}$ and the norm of its covariant derivative is at most $c_0 \frac{1}{\text{dist}(q, \cdot)^3}$. Multiply both sides of (6.13) by $\chi_{\diamond\diamond} \mathcal{G}_q$ with $\chi_{\diamond\diamond}$ being a smooth, non-negative function that is equal to 1 on $B_{r_\diamond}$ and equal to zero on the complement of the concentric ball with radius $\frac{33}{32} r_\diamond$. (Require also that $|d\chi_{\diamond\diamond}|$ be less than $c_0 \frac{1}{r_\diamond}$ and that $|\nabla d\chi_{\diamond\diamond}| \leq c_0 \frac{1}{r_\diamond^2}$.) Having done this multiplication, then integrate both sides of the resulting identity over the radius $\frac{33}{32} r_\diamond$ ball centered at p. Some integration by parts leads to an inequality that has the form

$$|v_n - v_\ddagger|(q) \leq c_0 \frac{1}{r_\diamond^4} \int_{B_{33r_\diamond/32}} |v_n - v_\ddagger|^2 + \int_{B_{33r_\diamond/32}} \frac{1}{\text{dist}(q, \cdot)^2} |\mathfrak{e}_n|^2 \ . \quad (6.15)$$

The right hand side of the various $n \in \Lambda$ versions of this inequality converge to zero as n gets ever larger; and (by virtue of (6.14) and the $L^2$ convergence of $\{v_n\}_{n \in \Lambda}$ to $v_\ddagger$), this convergence is uniform as q varies in $B_{r_\diamond}$.

*Part 4*: The strong $L^2_2$ convergence of $\{v_n\}_{n \in \Lambda}$ to $v_\ddagger$ follows from the $\{r_n, (A_n, a_n)\}_{n \in \Lambda}$ and $r = r_\diamond$ versions of the third bullet of Proposition 5.4. Proposition 5.4 can be invoked termwise for this sequence because of the assumption with regards to this part of Theorem 1.2 concerning the corresponding sequence of curvature $L^2$ norms.

**e) Verification of the assertions of Theorem 1.2**

This part of the subsection verifies the assertions of Theorem 1.2. To begin: Proposition 6.1 supplies a function $|v_\ddagger|$ and a subsequence $\Lambda \subset \mathbb{N}$ such that $\{|a_n|\}_{n \in \Lambda}$ converges to $|v_\ddagger|$ in the $L^2_1$ topology on compact subsets of X. Proposition 6.5 asserts that $\{|a_n|\}_{n \in \Lambda}$ also converges to $|v_\ddagger|$ in the $C^0$ topology on compact subsets of X.

Supposing henceforth that $|v_\ddagger|$ is not identically zero, the four parts that follow in this subsection verify the bulleted assertions of Theorem 1.2.

*Part 1*: Let Z denote the zero locus of $|v_\ddagger|$. The assertion made by the first bullet of Theorem 1.2 follows from the second and third bullets of Theorem 1.2 via an appeal to the theorems in [T1]. (The second and third bullets say in effect that the data set $(Z, \mathcal{I}, v_\ddagger)$ defines what [T1] calls a $\mathbb{Z}/2$ harmonic spinor. Meanwhile, theorems in [T1] say that the zero locus of a $\mathbb{Z}/2$ harmonic spinor is nowhere dense and has Hausdorff dimension at most 2. (These theorems say somewhat more about the zero locus.)

*Part 2*: With regards to the second bullet, Proposition 6.2 proves that $|v_\ddagger|$ is uniformly Hölder continuous on Z in compact subsets of X. If p is not in Z, then by virtue of what is said in Proposition 6.6, there is a ball in X−Z centered at p where $|v_\ddagger|$ is the norm of harmonic, self-dual 2 form. Since a harmonic differential form is $C^\infty$, its norm where it isn't zero is $C^\infty$ and therefore, $|v_\ddagger|$ is also $C^\infty$.

*Part 3*: This part of the proof proves the third bullet of Theorem 1.2. To this end, fix $\mu \in (0, \frac{1}{100}]$ and then some very large number for c. It should be greater than $c_\mu$ so that the propositions in Sections 4 and 5 and Lemma 6.4 and Proposition 6.6 can be used. (When X is non-compact, the open set $\underline{X}$ must also be fixed in advance.) Note that if $c > c_0$, then the mutual intersection of any two or three balls with centers in X ($\underline{X}$ if X is non-compact) will be convex. It simplifies matters to take c so that this is the case.

Given $p \in X-Z$ (in $\underline{X}-Z$ if X is not compact), let $r_p$ denote half of the smaller of $\liminf_{n \in \Lambda} r_{c\mu, n}$ and $\frac{1}{c}$. Lemma 6.4 says that this number is positive. The collection of balls labeled by the points in X−Z (or $\underline{X}$−Z as the case may be) with radius equal to the corresponding $r_p$ define an open cover of X−Z (or $\underline{X}$−Z). Fix a locally finite subcover to be denoted by $\mathfrak{U}$. Proposition 6.6 can be invoked on any given ball from $\mathfrak{U}$. Letting $B \subset \mathfrak{U}$ denote such a ball, let $v_{\ddagger B}$ denote the closed self-dual 2-form that is supplied by this proposition. This 2-form $v_{\ddagger B}$ is nowhere zero on B since its norm is $|v_\ddagger|$.

Let B and B´ denote intersecting balls from $\mathfrak{U}$. It follows from the fourth bullet of Propostion 6.1 and from the fourth bullet of Proposition 5.2 that $v_{\ddagger B} = \pm v_{\ddagger B'}$ on $B \cap B'$.

This + or - is the same at each point in B∩B´ because the latter set is connected. Use $\varepsilon_{BB'}$ to denote +1 if the + sign appears, or -1 if the - sign appears. The collection $\{\varepsilon_{BB'}\}_{B,B'\in\mathfrak{U}}$ with the open cover $\mathfrak{U}$ defines the cocycle data for a real line bundle over X–Z. This real line bundle is the bundle $\mathcal{I}$ of Theorem 1.2.

The data $\{v_{\ddagger B}\}_{B\in\mathfrak{U}}$ is the cocycle data for a section of $\Lambda^+\otimes \mathcal{I}$ (this follows tautologically from the definition of $\mathcal{I}$). This section of $\Lambda^+\otimes \mathcal{I}$ is what Theorem 1.2 denotes by $v_\ddagger$. It is a nowhere zero section because none of the $v_{\ddagger B}$'s are zero; and it is closed because $dv_{\ddagger B} = 0$ for each B.

*Part 4*: The assertions in Item a) of the fourth bullet of Theorem 1.2 follow from the versions of Proposition 5.2 and (5.3) that are obtained by using balls from $\mathfrak{U}$ and the data $(r = r_n, (A = A_n, a = a_n)$. The assertions in Item b) follow from the same versions of Proposition 5.4.

**7. Generalizations of (1.1) and (1.2)**

This section very briefly discusses the generalizations of (1.1) that are depicted in (1.4) and (1.6). Section 7a considers the former and Section 7b considers the latter.

**a) The equations in (1.4)**

As noted in Section 1's remark about (1.4), the appearance of $m\hat{a}$ in (1.4) (or, equivalently $rma$ in (1.5)) does not change the conclusions of either Theorem 1.1 or Theorem 1.2. To explain why, suppose that $r > 0$ and that $(A, a)$ obeys (1.5). The Bochner-Weitzenboch formula writes the equation $d_A^\dagger d_A a = 0$ as:

$$\nabla_A^\dagger \nabla_A a + r^2[a_c,[a,a_c]] + 2mr[a;a] = 0 .$$
(7.1)

This is the analog of (2.1). Taking the inner product of this equation with $a$ gives

$$\tfrac{1}{2}d^\dagger d|a|^2 + |\nabla_A a|^2 + r^2|[a;a]|^2 + 2mr\langle a,[a;a]\rangle + \langle a, \mathfrak{R}a\rangle = 0;$$
(7.2)

which is the analog of (2.2). Therefore, supposing that $\delta > 0$ is given, this identity (with the triangle inequality) leads directly to two inequalities:

- $\tfrac{1}{2}d^\dagger d|a|^2 + |\nabla_A a|^2 + (1-\delta)r^2|[a;a]|^2 \leq (\tfrac{m^2}{\delta} + |\mathfrak{R}|)|a|^2$ .
- $\tfrac{1}{2}d^\dagger d|a|^2 + |\nabla_A a|^2 + (1+\delta)r^2|[a;a]|^2 \geq -(\tfrac{m^2}{\delta} + |\mathfrak{R}|)|a|^2$ .

(7.3)

What with these two inequalities (in the case $\delta = \frac{1}{2}$ for example), all of the analysis that proves Theorems 1.1 and 1.2 can be repeated with almost no changes. This is literally true for the proof in [Ma] of Theorem 1.1. (Mare's proof is sketched in Section 2b). Meanwhile, the only change for the proof of Theorem 1.2 that is not completely straightforward is with regards to the definitions of N and K in Section 3a. With no modification of $\eth$ in (3.2), then the formula in (3.4) for N is different; it is:

$$N(r) = \frac{1}{r^2 K(r)^2} \int_{B_r} (|\nabla_A a|^2 + r^2 \, |[a;a]|^2 + 2mr\langle a, [a;a]\rangle)$$

(7.4)

This version of N need not be strictly positive. Even so, $N(r) \geq -c_0 r^2$ (because of (7.3)). All of the arguments in the proof of Theorem 1.2 go through with only minor cosmetic changes if N is allowed this small negativity. One can also make non-negative version of N by adding an $\mathcal{O}(r^2 \frac{m^2}{\delta} K^2)$ term to the definition in (3.2) of $\eth$. Indeed, by virtue of (7.3), the addition of a term

$$z \, \frac{m^2}{\delta} \int_0^r \left(\frac{1}{h(s)} \left(\int_{B_s} |a|^2\right)\right) ds$$

(7.5)

with $z \sim c_0$ will result in a non-negative version of N.

**b) The equations in (1.6)**

Suppose that $\hat{X}$ is a non-compact, oriented, Riemannian 4-manifold and suppose that $(A, (\hat{a}, \hat{\phi}))$ is a solution to (1.6) on $\hat{X}$. Fix an open set in $\hat{X}$ with compact closure to play the role of the manifold X. Because X has compact closure, the integrals of the square of the norm of both $\hat{a}$ and $\hat{\phi}$ on X are guaranteed to be finite. This understood, let $r$ denote the positive number whose square is given by

$$r^2 = \int_X (|\hat{a}|^2 + |\hat{\phi}|^2) \, .$$

(7.6)

Set $a = r^{-1}\hat{a}$ and set $\phi = r^{-1}\hat{\phi}$. The equations in (1.6) when written using $(a, \phi)$ read:

- $d_A a + *d_A \phi = 0.$
- $F_A^+ = \frac{1}{2} r^2 [a;a] + \frac{1}{\sqrt{2}} r^2 [\phi, a] \, .$
- $\int_X (|a|^2 + |\phi|^2) = 1 \, .$

(7.7)

Introduce by way of notation $q$ to denote the pair $(a, \phi)$, a section of $(\Lambda^+ \oplus \underline{\mathbb{R}}) \otimes \mathfrak{G}$. (The product real line bundle is denoted by $\underline{\mathbb{R}}$.) The equations in (7.7) are viewed as equations for a data set consisting of the positive number $r$ and the pair $(A, q)$.

The analog of Theorems 1.1 and 1.2 for the equations in (7.7) follows:

**Theorem 7.1**: *Let $\{(r_n, (A_n, a_n, \phi_n)\}_{n \in \mathbb{N}}$ denote a sequence of solutions to (7.7) on X.*

A. THE CASE $\{r_n\}_{n \in \mathbb{N}}$ HAS A CONVERGENT SUBSEQUENCE: *There is a principle SU(2) or SO(3) bundle $P' \to X$ (as the case may be), a non-negative number $r'$, a pair $(A', q')$ of solutions to (7.7) for the bundle $P'$, and a finite set of points in X (to be denoted by $\Xi$). These have the following significance: There is a subsequence $\Lambda \subset \mathbb{N}$ such that $\{r_n\}_{n \in \Lambda}$ converges to $r'$. And, there is a sequence of bundle isomorphisms $\{g_n : P'|_{X-\Xi} \to P|_{X-\Xi}\}_{n \in \Lambda}$ such that $\{(g_n{}^*A_n, g_n{}^*q_n\}_{n \in \Lambda}$ converges to $(A', q')$ in the $C^\infty$ topology on compact subsets of $X - \Xi$.*

B. THE CASE WHERE $\{r_n\}_{n \in \mathbb{N}}$ HAS NO CONVERGENT SUBSEQUENCE: *There exists a continuous and $L^2_{1;loc}$ function on X to be denoted by $|s_\ddagger|$ and a subsequence $\Lambda \subset \mathbb{N}$ such that $\{|q_n|\}_{n \in \Lambda}$ converges to $|s_\ddagger|$ in the $L^2_1$ and $C^0$ topologies on compact subsets of X. If $|s_\ddagger|$ is not identically zero, then the what follows is true.*
- *The zero locus of $|s_\ddagger|$ (which is denoted henceforth by Z) is a nowhere dense set with Hausdorff dimension at most 2.*
- *The function $|s_\ddagger|$ is smooth on the complement of Z and uniformly Hölder continuous along Z in compact subsets of X.*
- *There is a real line bundle $\mathcal{I} \to X-Z$ and a smooth section of $(\Lambda^+ \oplus \underline{\mathbb{R}}) \otimes \mathcal{I}$ on $X-Z$ whose norm is $|s_\ddagger|$. This section is harmonic in the following sense: Letting $s_\ddagger$ denote the section, write its components as $(\nu_\ddagger, \varphi_\ddagger)$. Then $d\nu_\ddagger + *d\varphi_\ddagger = 0$.*
- *There is a $\Lambda$-indexed sequence of isometric homomorphisms from $\mathcal{I}$ to $\mathfrak{G}|_{X-Z}$ to be denoted by $\{\sigma_n\}_{n \in \Lambda}$ that obeys the following:*
    a) *The sequence $\{q_n - s_\ddagger \sigma_n\}_{n \in \Lambda}$ converges to zero in the $C^0$ topology on compact subsets of $X-Z$. In addition, if $U \subset X-Z$ is any open set with compact closure, then*

    i) $\lim_{n \in \Lambda, n \to \infty} \int_U |\nabla_{A_n} \sigma_n|^2 = 0.$

    ii) $\lim_{n \in \Lambda, n \to \infty} \int_U |\nabla_{A_n}(q_n - s_\ddagger \sigma_n)|^2 = 0 .$

iii) $\lim_{n \in \Lambda, n \to \infty} \int_U |F_{A_n}^+|^2 = 0$.

b) *If there is an upper bound for the sequence* $\{\frac{1}{r_n^2} \int_X |F_{A_n}^-|^2\}_{n \in \Lambda}$ *and supposing that* $U \subset X-Z$ *is an open set with compact closure, then*

i) $\lim_{n \in \Lambda, n \to \infty} \int_U |\nabla_{A_n} \nabla_{A_n} \sigma_n|^2 = 0$.

ii) $\lim_{n \in \Lambda, n \to \infty} \int_U |\nabla_{A_n} \nabla_{A_n} (q_n - s_{\ddagger} \sigma_n)|^2 = 0$.

Part A of Theorem 7.1 can be proved by copying almost verbatim Mare's proof in [Ma] of Theorem 1.1. The only change (almost) is to replace $a$ everywhere by $q$. Note in this regard that the key Bochner-Weitzenboch formula for solutions to (7.7) (the analog of formula in (2.1) for solutions to (1.2)) has $\Lambda^+ \otimes \mathfrak{G}$ and $\mathfrak{G}$ components that read:

- $\nabla_A^\dagger \nabla_A a + r^2 [a_c, [a, a_c]] + r^2 [\varphi, [a, \varphi]] + \mathfrak{R} a = 0$.
- $\nabla_A^\dagger \nabla_A \varphi + r^2 [a_c, [\varphi, a_c]] = 0$.

(7.8)

These equations then lead to analogs of (2.2):

- $\frac{1}{2} d^\dagger d |a|^2 + |\nabla_A a|^2 + r^2 |[a; a]|^2 + r^2 |[\varphi, a]|^2 + \langle a, \mathfrak{R} a \rangle = 0$.
- $\frac{1}{2} d^\dagger d |\varphi|^2 + |\nabla_A \varphi|^2 + r^2 |[\varphi; a]|^2 = 0$.

(7.9)

(Note that integrating the lower equation on a compact X leads to the conclusion noted in Section 1 that both $\nabla_A \varphi$ and $[a, \varphi]$ must be zero.)

Part B of Theorem 7.1 can be proved by copying almost verbatim the proof of Theorem 1.1 in Sections 2-6 of this paper. One need only change $a$ to $q$ everywhere (and where $q$ is written using an orthonormal frame for $\Lambda^+$, one must write $q$ using an orthonormal frame for $\Lambda^+ \oplus \mathbb{R}$.) Some of the critical instances are noted below.

SECTION 3: With regards to the functions κ and N from Section 3a: The value of the relevant version of the function $\kappa^2$ for a given point $p \in X$ and value of r is (up to an $\mathcal{O}((1 + c_0 r^2))$ factor), the average of $\frac{1}{2\pi^2} |q|^2$ on the radius r sphere centered at p. Meanwhile, N is defined using this new version of κ; it is:

$$N(r) = \frac{1}{r^2 K(r)^2} \int_{B_r} (|\nabla_A a|^2 + r^2 \|[a;a]\|^2 + 2r^2 \|[a,\varphi]\|^2) \ .$$

(7.10)

With regards to the components of $q$: Suppose that $p \in X$ and that $e$ is a unit length element in the fiber of $\Lambda^+ \oplus \mathbb{R}$ at p. Parallel transport of $e$ on the radial geodesic arcs from p using the Levi-Civita connection defines a unit length section of $\Lambda^+ \oplus \mathbb{R}$ on the radius $r_0$ ball centered at p. Let $q_e$ denote the section of $\mathfrak{G}$ on this ball given by the inner product pairing between $e$ and $q$. This section $q_e$ obeys (by virtue of (7.8)) an equation that has the schematic form

$$\nabla_A^\dagger \nabla_A q_e + r^2 [q_c, [q_e, q_c]] + \Gamma_e \cdot \nabla_A q + \Upsilon_e \cdot q = 0 \ .$$

(7.11)

The notation here has $\{q_c\}_{c=1,2,3,4}$ being defined by the rules $\{q_c = a_c\}_{c=1,2,3}$ and $q_4 = \varphi$; and the convention is that repeated indices (the index c in this case) are summed. Meanwhile, $\Gamma_e$ and $\Upsilon_e$ are tensors that are defined from the Riemannian metric. They are $q$-independent and they obey $|\Upsilon_e| \le c_0$ and $|\Gamma_e| \le c_0 r$ on the ball $B_r$ if $r \in (0, r_0]$. This equation is the analog of (3.18) for solutions to (7.7).

SECTION 4: Fix $\mu \in (0, \frac{1}{100}]$ and $c > 1$ and a point $p \in X$. The number $r_{c\mu}$ is defined from the data $(r, (A, q))$ by the analog of (4.1) that replaces $|\nabla_A a|^2$ with $|\nabla_A q|^2$ and replaces $\|[a;a]\|^2$ with $\|[a;a]\|^2 + 2\|[\varphi;a]\|^2$. The wordings of the $(A, q)$ analogs of Propositions 4.1 and 4.2 stay the same. There are, however, two subtle points to be made with regards to the proof of Proposition 4.2. The first concerns an analog of (4.5) and thus (4.6) (which plays an absolutely central role). The analog is this: Let $e$ again denote a unit length section of $\Lambda^+ \oplus \mathbb{R}$ at a given point $p \in X$. Define, as before, $q_e$ to denote the metric induced pairing between $e$ and $q$. Let $q^\perp = q - e q_e$. (This is the part of $q$ that is annihilated by the pairing with $e$.) The analog of (4.5) has the schematic form

$$\nabla_A q_e + \mathcal{L}_e \nabla_A q^\perp = 0 \ ,$$

(7.12)

with $\mathcal{L}_e$ denoting a homomorphism that is defined canonically using the metric. This form for the equation $d_A a + *d_A \varphi$ is a consequence of the following fact: The fibers of $\Lambda^+ \oplus \mathbb{R}$ can be viewed as the quaternions (with the imaginary quaternions being the $\Lambda^+$ factor), and likewise so can the fibers of $T^*X$. Moreover, quaternionic multiplication (from the right) on the two copies of the quaternions is intertwined by the symbol of the operator $*(d + *d) : C^\infty(X; \Lambda^+ \oplus \mathbb{R}) \to C^\infty(X; T^*X)$.

The second subtle point concerns the identity $\frac{1}{4}|[a;a]|^2 = |a|^4 - \langle a_a a_b \rangle \langle a_a a_b \rangle$ that is used to deduce (4.23) from (4.22). There is a $q$ analog of this identity:

$$|q|^4 - \langle q_a q_b \rangle \langle q_a q_b \rangle = \frac{1}{4}(|[a;a]|^2 + 2|[\varphi,a]|^2).$$
(7.13)

(Note with regards to the derivation that $\langle [a;a],[\varphi,a] \rangle$ is identically zero.)

Granted (7.12) and (7.13), then the proof of the analog of Proposition 4.2 for solutions to (7.7) copies the arguments in Section 4 for the latter's version of Proposition 4.2. (Make sure however, to define the matrix $\mathbb{T}$ from (4.9) using all of the components of $q$. Thus, $\mathbb{T}$ should be a $4 \times 4$ matrix whose components, when written using an orthonormal frame for $\Lambda^+ \oplus \mathbb{R}$ at the given point, are defined by using $\{\langle q_a q_b \rangle\}_{a,b=1,2,3,4}$ in (4.9) in lieu of $\{\langle a_a, a_b \rangle\}_{a,b=1,2,3}$.)

SECTION 5: The analog for solutions to (7.7) of Proposition 5.1 is obtained from Proposition 5.1 by replacing $\{a_c\}_{c=1,2,3}$ in the statement by $\{q_c\}_{c=1,2,3,4}$. By virtue of (7.11), the latter can be the components of $q$ with respect to any orthonormal frame for $\Lambda^+ \otimes \mathbb{R}$ on the ball in question that is obtained by parallel transport on radial geodesics from an orthonormal frame at the ball's center point.

The analog of Proposition 5.2 for solutions to (7.7) writes $q$ as $q = s\sigma_{\ddagger} + \mathfrak{q}$ with $s$ being a section of $\Lambda^+ \oplus \mathbb{R}$ over the ball in question, and with $\mathfrak{q}$ obeying $\langle \sigma_{\ddagger} \mathfrak{q} \rangle = 0$ and $\langle s, \mathfrak{q} \rangle = 0$ on the ball. The analog of the second bullet of Proposition 5.2 is the same as in the original with $\mathfrak{a}$ replaced by $\mathfrak{q}$. The analog of Proposition 5.2's third bullet asserts that

$$\mathcal{D}s = \langle \mathfrak{Q}(\nabla_A \sigma_{\ddagger}) \mathfrak{q} \rangle$$
(7.14)

where $\mathcal{D}$ is the operator $d + *d$ (mapping $C^\infty(X; \Lambda^+ \oplus \mathbb{R})$ to $C^\infty(X; \wedge^3 T^*X)$), and where $\mathfrak{Q}$ is its principle symbol. The analog of the fourth and fifth bullets of Proposition 5.2 are identical with those in the original but for $\mathfrak{q}$ replacing $\mathfrak{a}$. Likewise, the analog of the sixth bullet is identical to the original with $s$ replacing $v$. The proof of this analog of Proposition 5.2 copies the proof of Proposition 5.2 with replacements that are straightforward and are left to the reader with no further comments.

There are likewise straightforward analogs of Propositions 5.4 and 5.5.

SECTION 6: All of the modifications to this section are straightforward replacements of $a$ by $q$, and $v$ by $s$, and $\mathfrak{a}$ by $\mathfrak{q}$; and nothing more will be said.

# References

[AF]  R. A. Adams and J. J. F. Fournier, Sobolev Spaces, Academic Press 2003

[Al]  F. J. Almgren, Jr. *Dirichlet's problem for multiple valued functions and the regularity of mass minimizing integral currents*, in Minimal Submanifols and Geodesics, Procedings of the Japan-United States Seminar, Tokyo 1977, pages 1-6. North Holland 1979.

[Ar]  N. Aronszajn, *A unique continuation theorem for solutions of elliptic partial differential equations or inequalities of second order*, J. Math. Pures. Appl. **36** (1957) 235-249.

[BW]  J. Bryan and R. Wentworth, *The multi-monopole equations for Kähler surfaces*, Turkish J. Math. **20** (1996) 119-128.

[DeL]  C. De Lellis, *The size of singular sets of area-minimizing currents*. arXiv1506.08118.

[DF]  H. Donnelly and C. Fefferman, *Nodal sets of eigenfunctions on Riemannian manifolds*, Invent. Math. **93** (1988) 161-183.

[HHL]  Q. Han, R. Hardt and F. Lin, *Geometric measure of singular sets of elliptic equations*, Commun. Pure and Applied Math., **51** (1998) 1425-1443.

[HW]  A. Haydys and T. Walpuski, *A compactness theorem for the Seiberg-Witten equations with multiple spinors in dimension three.*, Geom. Funct. Anal. **25** (2015) 1799-1821 (Erratum available via doi:10.13140/RG.2.1.2559.1285.

[Ma]  B. Mares, *Some analytic aspects of Vafa-Witten twisted N =4 supersymmetric Yang-Mills theory*, Ph. D. thesis, M.I.T. 2010.

[Mo]  C. B. Morrey, Multiple Integrals in the Calculus of Variations, Springer 1966.

[Ni]  L. Nicolescu, Lectures on the Geometry of Manifolds (2$^{nd}$ edition), World Scientific 2007.

[Ta]  Y. Tanaka, *Some boundedness properties of solutions to the Vafa-Witten equations on closed 4-manifolds*, arXiv.1308.0862.

[T1]  C. H. Taubes, *The zero loci of $\mathbb{Z}/2$ harmonic spinors in dimension 2, 3 and 4*. arXiv:1407.6206.

[T2]  C. H. Taubes, *Compactness theorems for SL(2; $\mathbb{C}$) generalizations of the 4-dimensional anti-self dual equations*. arXiv:1307.6447.

[T3]  C. H. Taubes, *On the behavior of sequences of solutions to U(1) Seiberg-Witten systems in dimension 4*. arXiv:1610.07163.

[T4]  C. H. Taubes, *PSL(2; $\mathbb{C}$) connections with $L^2$ bounds on curvature*, Cambr. J. Math **1**  (2013) 239-397 and *Corrigendum to "PSL(2; $\mathbb{C}$) connections with $L^2$ bounds on curvature"*, Cambr. J. Math **3** (2015) 619-631.

[T5]  C. H. Taubes, *Growth of the Higgs field for solutions to the Kapustin-Witten equations on $\mathbb{R}^4$*.  arXiv:1701.03072.

63
[U1]   K. K. Uhlenbeck, *Connections with $L^p$ bounds on curvature*, Commun. Math. Phys. **83** (1982) 31-42.

[U2]   K. K. Uhlenbeck, *Removable singularities in Yang-Mills fields*, Commun. Math. Phys. **83** (1982) 11-29.

[VW]   C. Vafa and E. Witten, *A strong coupling test of S-duality*, Nucl. Phys. B. 431 (1994) 3-77.